\documentclass[11pt]{article}
\usepackage{geometry}
\usepackage{amssymb}
\usepackage{xcolor}
\usepackage{amsmath}
\usepackage{mathrsfs}
\usepackage{amsthm}
\allowdisplaybreaks[4]

\begin{document}
\flushleft

\begin{center}
\Large\textbf{On a variance associated with the distribution of $r_3(n)$ in arithmetic progressions}

\end{center}

\begin{center}

\vspace{3mm}

\large
{Pengyong Ding

Department of Mathematics, Pennsylvania State University}

\end{center}

\begin{center}

\large\textbf {1. Introduction}

\end{center}

This paper is concerned with the conjunction of two questions which have attracted much attention over the years.  The first question, having its origin in work of Barban \cite{Bar} in the special case $a_n=\Lambda(n)$, the von Mangoldt function and largely stimulated by developments in the large sieve, concerns the variance

\begin{equation}\label{vxqoriginal}
V(x, Q)=\sum_{q\le Q}\sum^{q}_{a=1} |A(x;q,a)-f(q,a)M(x)|^2,
\end{equation}

where

\begin{equation}\label{axqaoriginal}
A(x;q,a)=\sum_{\substack{
n\le x\\
n\equiv a(\text{mod} q)
}} a_n,
\end{equation}

$\{a_n\}$ is a real sequence, and $f(q,a)$ and $M(x)$ approximately reflect the local and global properties of $\{a_n\}$ respectively.  The special case was studied successively by Davenport and Halberstam \cite{DH}, and Gallagher \cite{Ga}, Montgomery \cite{Mon} and Hooley \cite{Ho}.  Montgomery established an asymptotic formula in the special case and a more refined result with a simplified proof was given by Hooley.  Hooley \cite{Ho2}, \cite{Ho3}, \cite{Ho5}, \cite{Ho6}, \cite{Ho7}, \cite{Ho8}, \cite{Ho9} then developed the subject further by studying wide reaching generalizations.  In particular his method established that an asymptotic formula can be obtained for a wide class of $a_n$ provided that only that one has an asymptotic formula of the kind
\begin{equation}
\label{SiegelWalfisz}
A(x;q,a)\sim xf(q,a)
\end{equation}
with a reasonable error term, analogous to the Siegel-Walfisz theorem, and some control over the behaviour of
\begin{equation}
\label{l2mean}
\sum_{n\le x}a_n^2.
\end{equation}
In the special case of the von Mangoldt function further refinements occur in Friedlander and Goldston \cite{FG}, and in Vaughan and Goldston \cite{GV}.  The latter paper showing that some advantage could be accrued by applying a version of the Hardy-Littlewood method.  Having shown the usefulness of the this method in that specific question, Vaughan \cite{Va3},\cite{Va4} then considered the more general problem.  Suppose that
\[
M(x)=x
\] and $f(q,a)$ only depends on $q$ and $(q,a)$, and the real sequence $\{a_n\}$ satisfies \[
\sum_{n\le x} a_n^2\ll x.
\]
For convenience write $f(q, (q,a))$ for $f(q,a)$ and suppose further that
\[
A(x;q,a)=xf(q, (q,a))+ O\Big(\frac{x}{\Psi(x)}\Big),
\]
where $\Psi(x)$ is an increasing function with $\Psi(x)>\log x$ for all large $x$, $\Psi(1)>0$ and
\[
\int_{1}^{x} \Psi(y)^{-1} \mathrm{d}y\ll x\Psi(x)^{-1}
\]
Then it was shown that
\[
V(x,Q)=Q\sum_{n\le x}a_n^2-Qx\sum^{\infty}_{q=1} g(q)+U(x,Q),
\]
where
\[
U(x,Q)\ll x^{3/2}\log x+x^2(\log 2x)^{9/2} \Psi(x)^{-1}+x^2(\log x)^{4/3}\Psi(x)^{-2/3}+Q^2 E(x/Q),
\]
\[
E(z)=\displaystyle\int_0^{z} \sum_{q>y}^{} g(q) \mathrm{d}y,
\]
and
\[g(q)=\displaystyle\phi(q)\Big(\sum_{r|q} f(q,r)\mu(q/r)\Big)^2.
\]
Note that the distribution of $\Lambda(n)$ is not included in this case, since
\[
\sum_{n\le x} \Lambda(n)^2=x\log x+O(x),
\]
contrary to the requirement that $\displaystyle\sum_{n\le x} a_n^2\ll x$.
\par
An example of an arithmetical function whose mean value is a constant but whose means square is unbounded is $r(n)$, the number of ways of writing a positive integer $n$ as the sum of two squares.  This was studied in the context of the above variance by Dancs \cite{Dan}.   It is, therefore of some interest to consider examples in which the mean square of $a_n$, and perhaps also the mean value are unbounded.  Thus with this in mind the above variance has been studied in the special case of $a_n=d(n)$, the divisior function, by Pongsriiam \cite{P} and Pongsriiam and Vaughan \cite{PV}.
\par
In all of the special cases studied hitherto, the behaviour of the mean square is well understood.  The second question which we conjoin to the first concerns the function \begin{equation}\label{defdefdefr3n}
r_3(n)=\displaystyle\sum_{\substack{
x_1, x_2, x_3 \\
x_1^3+x_2^3+x_3^3=n
}} 1,
\end{equation}
the number of ordered representations of $n$ as the sum of three positive cubes.  In spite of significant work by Hooley this function is only poorly understood.  We know (Vaughan \cite{RV20}) that

\begin{equation}\label{r3nsum}
\sum_{n\le x} r_3(n)= \Gamma\Big(\frac{4}{3}\Big)^3 x-\frac{\Gamma(\frac{4}{3})^2}{2\Gamma(\frac{5}{3})} x^{\frac{2}{3}} +O(x^{\frac{5}{9}} (\log x)^{\frac{1}{3}}),
\end{equation}

and 

\begin{equation}\label{r3nsquaredsum}
\sum_{n\le x} r_3(n)^2\ll x^{\frac{7}{6}} (\log x)^{\varepsilon-\frac{5}{2}},
\end{equation}
and assuming that a certain Hasse-Weil $L$-function satisfies the Riemann Hypothesis Hooley \cite{CH96} has shown that
\[
\sum_{n\le x} r_3(n)^2\ll x^{1+\varepsilon}.
\]
\par
It may be true that
\[
\sum_{n\le x} r_3(n)^2\sim Cx
\]
for some positive constant $C$, but we only know a lower bound of this order with a value of $C$ larger than the obvious guess.  See Hooley \cite{CH86}.  Also,  by refining the deep work of Hooley \cite{CH81}, Vaughan has shown that for every $\varepsilon>0$ there exists a positive constant $\delta$, such that for every $\varepsilon>0$ and all $Q\le x$,
\begin{equation}\label{estimateupsilonxqa}
\sum_{q\le Q} \max_a \sup_{y\le x} \left|
\Upsilon(y;q,a)-\Gamma\Big(\frac{4}{3}\Big)^3y\rho(q,a)q^{-3}
\right|\ll_{\varepsilon} x^{\frac{8}{9}+\varepsilon}+x^{\frac{1}{3}}Q^{\frac{2}{9}}(Q^{\frac{10}{9}}+x^{\frac{5}{9}})(\log x)^{-\delta},
\end{equation}

where

\begin{equation}\label{upsilonxqadefinition}
\Upsilon(x;q,a)=\sum_{\substack{
n\le x\\
n\equiv a(\text{mod} q)
}} r_3(n).
\end{equation}

and $\rho(q,a)$ denotes the number of solutions of the congruence $l_1^3+l_2^3+l_3^3\equiv a (\text{mod} q)$.

\vspace{3mm}

When we compare (\ref{estimateupsilonxqa}), (\ref{upsilonxqadefinition}) with (\ref{axqaoriginal}), and use the global property (\ref{r3nsum}) of $r_3(n)$, we have:

\begin{equation}\label{axqaforr3n}
A(x;q,a)=\Upsilon(x;q,a),
\end{equation}

\begin{equation}\label{fqaforr3n}
f(q,a)= \rho(q,a)q^{-3},
\end{equation}

\begin{equation}\label{mxforr3n}
M(x)= \Gamma\Big(\frac{4}{3}\Big)^3 x,
\end{equation}

and more the variance (\ref{vxqoriginal}) becomes

\begin{equation}\label{interestingvxq}
V(x, Q)=\sum_{q\le Q}\sum^{q}_{a=1} \left|
\Upsilon(x;q,a)-\Gamma\Big(\frac{4}{3}\Big)^3 x\rho(q,a)q^{-3}
\right|^2,
\end{equation}

when $a_n=r_3(n)$.

\vspace{3mm}

The purpose of this article is to use the Hardy - Littlewood Method and the Farey dissection to estimate $V(x,Q)$ stated above when $x$ is large enough and $1\le Q\le x$, especially when $Q$ is not too small (for example, $Q>\sqrt{x}$ or $Q>\sqrt{x}\log 2x$), otherwise the result is trivial. We will start from the results of the general problem (\ref{vxqoriginal}), after that, we will study the variance (\ref{interestingvxq}). For future reference, in this article we always assume that $x$ is sufficiently large.

\vspace{3mm}

\begin{center}

\large\textbf{2. The standard initial procedure}

\end{center}

From (\ref{vxqoriginal}) and (\ref{axqaoriginal}), we have

\begin{equation}\label{vxqoriginalfirstchange}
V(x,Q)=2S_1-2S_2+S_3+[Q] \sum_{n\le x} a_n^2,
\end{equation}

where

\begin{equation}\label{s1original}
S_1=\sum_{q\le Q}\sum_{\substack{
m<n\le x\\
m\equiv n(\text{mod} q)
}} a_ma_n,
\end{equation}

\begin{equation}\label{s2original}
S_2=M(x)\sum_{q\le Q}\sum_{a=1}^q f(q,a)\sum_{\substack{
n\le x\\
n\equiv a(\text{mod} q)
}} a_n,
\end{equation}

\begin{equation}\label{s3original}
S_3=M(x)^2\sum_{q\le Q}\sum_{a=1}^q f(q,a)^2.
\end{equation}

It is possible to rewrite $-2S_2+S_3$ as $-S_3+2(S_3-S_2)$, where $S_3-S_2$ is usually an error term. Also, in order to avoid $[Q]$ when $Q$ is not an integer, we would rather use $\displaystyle Q\sum_{n\le x} a_n^2+O\Big(\sum_{n\le x} a_n^2\Big)$ than $\displaystyle[Q]\sum_{n\le x} a_n^2$. Later, we will show that after using the Farey dissection, $\displaystyle O\Big(\sum_{n\le x} a_n^2\Big)$ will be dominated by the other error terms.

\vspace{3mm}

The most difficult part is the treatment of $S_1$. Let

\begin{equation}\label{fofalpha}
F(\alpha)=\sum_{q\le Q}\sum_{r\le x/q} e(\alpha qr),
\end{equation}

\begin{equation}\label{gofalpha}
G(\alpha)=\sum_{n\le x} a_ne(n\alpha),
\end{equation}

then

\begin{equation}\label{s1withfofalphaandgofalpha}
S_1=\int_0^1 F(\alpha)|G(\alpha)|^2\mathrm{d}\alpha.
\end{equation}

From the proof in Vaughan \cite{Va4}, we know that $F(\alpha)=F_q(\alpha)+H_q(\alpha)$, where

\begin{equation}\label{fsubqofalpha}
F_q(\alpha)=\sum_{\substack{
l\le\sqrt{x}\\
q|l
}}\Big(\sum_{m\le x/l}+\sum_{\sqrt{x}<m\le \min (Q, x/l)}\Big)e(\alpha lm)
\end{equation}

and $H_q(\alpha)$ is the corresponding multiple sum with $q\nmid l$.

\vspace{3mm}

Now we suppose that $R$ satisfies

\begin{equation}\label{R}
2\sqrt{x}\le R\le \frac{1}{2}x
\end{equation}

and consider the intervals $\mathfrak{M}(q,a)$ associated with the element $a/q$ of the Farey dissection of order $R$ when $1\le a\le q\le R$ and $(a,q)=1$. Then these intervals are pairwise disjoint, and the union of them is $(1/(R+1), (R+2)/(R+1)]$, which is an interval of length 1. By (\ref{s1withfofalphaandgofalpha}), we have

\begin{equation}\label{s1withfsubqofalphaandhsubqofalpha}
S_1=\sum_{q\le R}\sum_{\substack{
a=1\\ (a,q)=1
}}^q\int_{\mathfrak{M}(q,a)}F_q(\alpha)|G(\alpha)|^2\mathrm{d}\alpha+\sum_{q\le R}\sum_{\substack{a=1\\ (a,q)=1
}}^q\int_{\mathfrak{M}(q,a)}H_q(\alpha)|G(\alpha)|^2\mathrm{d}\alpha.
\end{equation}

The second term above is $\displaystyle\ll R(\log x) \sum_{n\le x} a_n^2$, since $H_q(\alpha)\ll R\log x$ when $\alpha\in\mathfrak{M}(q,a)$. So we can write it as $\displaystyle O\Big(R(\log x)\sum_{n\le x} a_n^2\Big)$, which dominates $O\Big(\displaystyle\sum_{n\le x} a_n^2\Big)$.

\vspace{3mm}

To estimate the first term, let $\beta=\alpha-a/q$ when $\alpha\in\mathfrak{M}(q,a)$, then

\begin{equation}\label{estimatefsubqofalpha}
F_q(\alpha)\ll \frac{x\log(2\sqrt{x}/q)}{q+qx|\beta|}
\end{equation}

if $q\le\sqrt{x}$ and $|\beta|\le \displaystyle\frac{1}{2q\sqrt{x}}$ from the proof in Goldston and Vaughan \cite{GV}.

\vspace{3mm}

Now define the major arc $\mathfrak{N}(q,a)$ by

\begin{equation}\label{definitionofnqa}
\mathfrak{N}(q,a)=\Big[\frac{a}{q}-\frac{1}{2qR},\frac{a}{q}+\frac{1}{2qR}\Big],
\end{equation}

then $\mathfrak{N}(q,a)\subset\mathfrak{M}(q,a)$. We can prove that $F_q(\alpha)\ll R\log x$ if $\alpha\in\mathfrak{M}(q,a)\setminus\mathfrak{N}(q,a)$, or if $\alpha\in\mathfrak{N}(q,a)$ but $q>x/R$. So

\begin{equation}\label{s1withs4}
S_1=S_4+O\Big(R(\log x) \sum_{n\le x} a_n^2\Big)
\end{equation}

where

\begin{equation}\label{definitionofs4}
S_4=\sum_{q\le x/R}\sum_{\substack{
a=1\\ (a,q)=1
}}^q\int_{\mathfrak{N}(q,a)}F_q(\alpha)|G(\alpha)|^2\mathrm{d}\alpha=\sum_{q\le x/R}\int_{-1/2qR}^{1/2qR}F_q(\beta)\sum_{\substack{a=1\\
(a,q)=1}}^q \Big|G\Big(\beta+\frac{a}{q}\Big)\Big|^2\mathrm{d}\beta.
\end{equation}

The last part of the equality comes from integration by substitution $\alpha=\beta+a/q$, and the property that $F_q(\alpha)$ is periodic with period $1/q$.

\vspace{3mm}

From (\ref{vxqoriginalfirstchange}) and (\ref{s1withs4}), we have

\begin{equation}\label{vxqoriginalwiths2s3s4}
V(x,Q)=2S_4-S_3+2(S_3-S_2)+Q \sum_{n\le x} a_n^2+O\Big(R\log x\sum_{n\le x}a_n^2\Big),
\end{equation}

where $\displaystyle Q\sum_{n\le x} a_n^2$ is expected to be one of the main terms at least when $Q/(R\log x)$ is large.  From now on, we always assume that
\[
Q\ge x^{1/2}\log x.
\]

\vspace{3mm}

The treatment of $S_4$ depends on obtaining an approximation to $G(\beta+a/q)$ whicn can be considered a product of local factors.  In the examples discussed above, it is possible to estimate $G(\beta+a/q)$ as $\nu(q)J(\beta)$ for some $\nu(q)$ and $J(\beta)$ that are independent of $a$, at least when $(q,a)=1$.  In the situation which arises with $r_3$ we have to deal with a more complex approximation $\nu(q,a)J(\beta)$.   Nevertheless, for a given $q$, it is possible to factor out
\[
\sum^q_{\substack{a=1\\
(a,q)=1
}} |\nu(q,a)|^2,
\]
and the integral will be changed to
\[
\int_{-1/2qR}^{1/2qR}F_q(\beta)|J(\beta)|^2\mathrm{d}\beta,
\]
which is likely to be close to
\[
\int_{-1/2}^{1/2} F_q(\beta)|J(\beta)|^2 \mathrm{d}\beta.
\]

\vspace{3mm}

\begin{center}

\large\textbf{3. Preliminary results of the variance (\ref{interestingvxq})}

\end{center}

If $a_n=r_3(n)$, then by \S2, the variance (\ref{interestingvxq}) is:

\begin{equation}\label{vxqwiths2s3s4forr3n}
V(x,Q)=2S_4-S_3+2(S_3-S_2)+Q \sum_{n\le x} r_3(n)^2+O\Big(R(\log x)\sum_{n\le x}r_3(n)^2\Big),
\end{equation}

where $S_4$ is still defined by (\ref{definitionofs4}) and $G(\alpha)$ has the expression:

\begin{equation}\label{gofalphaforr3n}
G(\alpha)=\sum_{n\le x} r_3(n)e(n\alpha).
\end{equation}

Also, by (\ref{fqaforr3n}), (\ref{mxforr3n}), (\ref{s2original}) and (\ref{s3original}), now $S_2$ and $S_3$  have the expressions:

\begin{equation}\label{s2forr3n}
S_2=\Gamma\Big(\frac{4}{3}\Big)^3 x\sum_{q\le Q}\sum_{a=1}^q \frac{\rho(q,a)}{q^3}\sum_{\substack{
n\le x\\
n\equiv a(\text{ mod } q)
}} r_3(n),
\end{equation}

\begin{equation}\label{s3forr3n}
S_3=\Gamma\Big(\frac{4}{3}\Big)^6 x^2\sum_{q\le Q}\sum_{a=1}^q \frac{\rho(q,a)^2}{q^6}.
\end{equation}

In the example of the distribution of $a_n=\Lambda(n)$, it is easy to get a good approximation for $\displaystyle\sum_{n\le x} a_n^2$, which is $\displaystyle\sum_{n\le x} \Lambda(n)^2=x\log x+O(x)$. However, when $a_n=r_3(n)$ and $n\ge 3$, our understanding of $r_3(n)$ is so poor that (\ref{r3nsquaredsum}) is all we know about $\displaystyle\sum_{n\le x} r_3(n)^2$. We also have $x\ll \displaystyle\sum_{n\le x} r_3(n)^2$ from (\ref{r3nsum}), since $r_3(n)$ is a nonnegative integer. First, we can redefine the upper exponent $A^+$ and lower exponent $A^-$ of $x$ in $\displaystyle\sum_{n\le x} r_3(n)^2$ as follows:

\begin{equation}\label{approximationr3nsquaredsumpositive}
A^+=\limsup_{x\to\infty} \frac{\displaystyle\log\sum_{n\le x} r_3(n)^2}{\log x},
\end{equation}

\begin{equation}\label{approximationr3nsquaredsumnegative}
A^-=\liminf_{x\to\infty} \frac{\displaystyle\log\sum_{n\le x} r_3(n)^2}{\log x},
\end{equation}

and if $A^+=A^-$, then redefine the exponent $A$ as follows:

\begin{equation}\label{approximationr3nsquaredsum}
A=A^-=A^+.
\end{equation}

From the discussion above, we can infer that

\begin{equation}\label{exponentofr3nsquaresum}
1\le A^- \le A^+\le \frac{7}{6}.
\end{equation}

Later in this article, we will show that the expression of $R$, which is the order of the Farey dissection, depends on the value of $A^+$. So the range of $Q$ such that the conclusion is not trivial as well as the error terms also depend on $A^+$. Once the value of $A^+$ is found after some research on $\displaystyle\sum_{n\le x} r_3(n)^2$, the conclusion mentioned in this article will be finalized automatically.

\vspace{3mm}

\begin{center}

\large\textbf{4. Several lemmata}

\end{center}

Here we prove some useful results. First, define a function $S(q,a)$ as follows:

\begin{equation}\label{sofqb}
S(q,a)=\sum^q_{m=1} e\Big(\frac{am^3}{q}\Big).
\end{equation}

\textbf{Lemma 1.} \textit{Let $\rho(q,a)$ denotes the number of solutions of the congruence $l_1^3+l_2^3+l_3^3\equiv a (\text{mod} q)$. Then}

\begin{equation}\label{rhoofqa}
\rho(q,a)=\frac{1}{q}\sum_{b=1}^q e\Big(-\frac{ba}{q}\Big) S(q,b)^3,
\end{equation}

\textit{and}

\begin{equation}\label{sumrhoofqasquared}
\sum_{a=1}^q \rho(q,a)^2=\frac{1}{q} \sum_{b=1}^q |S(q,b)|^6=q^5\sum_{r|q} \frac{1}{r^6}\sum_{\substack{c=1\\ (c,r)=1
}}^r |S(r,c)|^6.
\end{equation}

\begin{proof}
The value of $\displaystyle\frac{1}{q}\sum_{b=1}^q e\Big(\frac{b(l_1^3+l_2^3+l_3^3-a)}{q}\Big)$ is 1 if $\displaystyle\frac{l_1^3+l_2^3+l_3^3-a}{q}$ is an integer, and 0 otherwise. In other words,

\begin{equation*}
\frac{1}{q}\sum_{b=1}^q e\Big(-\frac{ba}{q}\Big)e\Big(\frac{bl_1^3}{q}\Big)e\Big(\frac{bl_2^3}{q}\Big)e\Big(\frac{bl_3^3}{q}\Big)=1
\end{equation*}

iff $(l_1, l_2, l_3) (\text{mod} q)$ is a solution of the congruence $l_1^3+l_2^3+l_3^3\equiv a (\text{mod} q)$. So the number of all such solutions is
\begin{align*}
\rho(q,a)&=\sum_{l_1=1}^q\sum_{l_2=1}^q\sum_{l_3=1}^q\Big(\frac{1}{q}\sum_{b=1}^q e\Big(-\frac{ba}{q}\Big)e\Big(\frac{bl_1^3}{q}\Big)e\Big(\frac{bl_2^3}{q}\Big)e\Big(\frac{bl_3^3}{q}\Big)\Big)\\
&=\frac{1}{q}\sum_{b=1}^q e\Big(-\frac{ba}{q}\Big)\Big(\sum_{l_1=1}^q e\Big(\frac{bl_1^3}{q}\Big)\Big)\Big(\sum_{l_2=1}^q e\Big(\frac{bl_2^3}{q}\Big)\Big)\Big(\sum_{l_3=1}^q e\Big(\frac{bl_3^3}{q}\Big)\Big)\\
&=\frac{1}{q}\sum_{b=1}^q e\Big(-\frac{ba}{q}\Big) S(q,b)^3.
\end{align*}

We have $\rho(q,a)=\overline{\rho(q,a)}=|\rho(q,a)|$ since $\rho(q,a)$ is a nonnegative integer. So
\begin{align*}
\sum_{a=1}^q \rho(q,a)^2&=\sum_{a=1}^q\Big(\frac{1}{q}\sum_{b_1=1}^q e\Big(-\frac{b_1a}{q}\Big) S(q,b_1)^3\Big)\Big(\frac{1}{q}\sum_{b_2=1}^q e\Big(\frac{b_2a}{q}\Big) \overline{S(q,b_2)^3}\Big)\\
&=\frac{1}{q^2}\sum_{b_1=1}^q\sum_{b_2=1}^q S(q, b_1)^3\overline{S(q, b_2)^3}\sum_{a=1}^q e\Big(\frac{(b_2-b_1)a}{q}\Big)\\
&=\frac{1}{q^2}\sum_{1\le b_1=b_2\le q} S(q, b_1)^3\overline{S(q, b_2)^3} \cdot q\\
&=\frac{1}{q}\sum_{b=1}^q |S(q,b)|^6.
\end{align*}

To complete the proof, note that if $c/r$ is the simplified form of $b/q$, in other words, if $c/r=b/q$ and $(c,r)=1$, then $S(q,b)=(q/r)S(r,c)$. So

\begin{equation*}
\sum_{b=1}^q|S(q,b)|^6=\sum_{r|q}\sum^r_{\substack{c=1\\
(c,r)=1}} \Big|\frac{q}{r}S(r,c)\Big|^6=q^6\sum_{r|q}\frac{1}{r^6}\sum_{\substack{c=1\\ (c,r)=1
}}^r |S(r,c)|^6.
\end{equation*}
\end{proof}

\vspace{3mm}

In order to estimate $S_3$, it is convenient to define a function $T(r)$ as follows:

\begin{equation}\label{tofr}
T(r)=\frac{1}{r^7}\sum^r_{\substack{c=1\\
(c,r)=1
}}|S(r,c)|^6.
\end{equation}

\vspace{3mm}

\textbf{Lemma 2.} We have $T(r)\ll r^{-2}$.

\begin{proof}
We have that if $(q,a)=1$, then $S(q,a)\ll q^{1-1/k}$. Now $k=3$, so $|S(r,c)|\ll r^{2/3}$ if $(c,r)=1$. So

\begin{equation*}
T(r)\ll \frac{1}{r^7}\sum_{\substack{c=1\\
(c,r)=1
}}^r (r^{2/3})^6=\frac{\phi(r)}{r^3}\le \frac{1}{r^2}.
\end{equation*}
\end{proof}

\vspace{3mm}

\textbf{Lemma 3.} \textit{Suppose that $(q,a)=1$, then for every $\varepsilon>0$,}

\begin{equation}\label{lemma3}
\sum_{n\le x} r_3(n)e\Big(\frac{an}{q}\Big)=\Gamma\Big(\frac{4}{3}\Big)^3xq^{-3}S(q,a)^3+O\Big(x^{\frac{2}{3}}q^{\frac{1}{2}+\varepsilon}\Big).
\end{equation}

\begin{proof}
By (\ref{defdefdefr3n}),
\[
r_3(n)=\displaystyle\sum_{\substack{x_1, x_2, x_3 \\ x_1^3+x_2^3+x_3^3=n
}} 1.
\]
So
\begin{align*}
\sum_{n\le x} r_3(n)e\Big(\frac{an}{q}\Big)&=\sum_{n\le x} \sum_{\substack{
x_1, x_2, x_3 \\
x_1^3+x_2^3+x_3^3=n
}} e\Big(\frac{a(x_1^3+x_2^3+x_3^3)}{q}\Big)\\
&= \sum_{\substack{x_1, x_2\\
x_1^3+x_2^3\le x
}} e\Big(\frac{a(x_1^3+x_2^3)}{q}\Big) \sum_{x_3\le (x-x_1^3-x_2^3)^{1/3}} e\Big(\frac{ax_3^3}{q}\Big).
\end{align*}

Consider the innermost sum. We have that if $(a,q)=1$, $\alpha=a/q+\beta$ and $|\beta|\le (2kq)^{-1}n^{1/k-1}$, then for every $\varepsilon>0$,

\begin{equation}\label{theorem4.1of[4]}
\sum_{x\le n^{1/k}} e(\alpha x^k)-q^{-1}S(q,a)v_1(\alpha-a/q) \ll q^{1/2+\varepsilon},
\end{equation}

where $v_1(\beta)=\displaystyle\int_0^{n^{1/k}} e(\beta\gamma^k)\mathrm{d}\gamma$. Now $k=3$, $\alpha=a/q$ (so $\beta=0$ satisfies the condition), and substitute $n$ with $x-x_1^3-x_2^3$, so $v_1(\alpha-a/q)=v_1(0)=(x-x_1^3-x_2^3)^{1/3}$. Hence (\ref{theorem4.1of[4]}) becomes

\begin{equation*}
\sum_{x_3\le (x-x_1^3-x_2^3)^{1/3}} e\Big(\frac{ax_3^3}{q}\Big) -q^{-1}S(q,a) (x-x_1^3-x_2^3)^{1/3}\ll q^{1/2+\varepsilon}.
\end{equation*}

In other words, we can write the innermost sum as $q^{-1}S(q,a)(x-x_1^3-x_2^3)^{1/3}+O(q^{1/2+\varepsilon})$. So
\begin{align}\label{lemma3nextstep}
\sum_{n\le x} r_3(n)e\Big(\frac{an}{q}\Big)&=q^{-1}S(q,a)\sum_{\substack{
x_1, x_2\\
x_1^3+x_2^3\le x
}} e\Big(\frac{a(x_1^3+x_2^3)}{q}\Big) (x-x_1^3-x_2^3)^{1/3}\notag\\
&+O \Big(\sum_{\substack{x_1, x_2\\ 
x_1^3+x_2^3\le x
}}  q^{1/2+\varepsilon}\Big).
\end{align}

The second term above is

\begin{equation*}
\ll q^{1/2+\varepsilon}\sum_{\substack{x_1, x_2 \\
x_1^3+x_2^3\le x
}} 1 \le q^{1/2+\varepsilon} x^{2/3}.
\end{equation*}

For the first term, repeat this argument on the sum over $x_2$. The first term is
\begin{align*}
&=q^{-1}S(q,a)\sum_{x_1\le x^{1/3}} e\Big(\frac{ax_1^3}{q}\Big) \sum_{x_2\le (x-x_1^3)^{1/3}} e\Big(\frac{ax_2^3}{q}\Big) (x-x_1^3-x_2^3)^{1/3}\\
&=q^{-1}S(q,a)\sum_{x_1\le x^{1/3}} e\Big(\frac{ax_1^3}{q}\Big) \sum_{x_2\le (x-x_1^3)^{1/3}} e\Big(\frac{ax_2^3}{q}\Big) \int_{x_2}^{(x-x_1^3)^{1/3}} (x-x_1^3-y^3)^{-2/3} y^2\mathrm{d}y\\
&=q^{-1}S(q,a)\sum_{x_1\le x^{1/3}} e\Big(\frac{ax_1^3}{q}\Big) \int_{0}^{(x-x_1^3)^{1/3}} (x-x_1^3-y^3)^{-2/3}y^2 \Big(\sum_{x_2\le y} e\Big(\frac{ax_2^3}{q}\Big)\Big)\mathrm{d}y
\end{align*}

by changing the order of the sum and the integral. For the innermost sum in the parenthesis, we need to substitute $n$ with $y^3$ in (\ref{theorem4.1of[4]}):

\begin{equation*}
\sum_{x_2\le y} e\Big(\frac{ax_2^3}{q}\Big)-q^{-1}S(q,a)y \ll q^{1/2+\varepsilon},
\end{equation*}

in other words, the innermost sum is $q^{-1}S(q,a)y+O(q^{1/2+\varepsilon})$. Hence the first term of (\ref{lemma3nextstep}) is
\begin{align}\label{lemma3nextnextstep}
&q^{-2}S(q,a)^2\sum_{x_1\le x^{1/3}} e\Big(\frac{ax_1^3}{q}\Big) \int_{0}^{(x-x_1^3)^{1/3}} (x-x_1^3-y^3)^{-2/3}y^3 \mathrm{d}y\notag\\
+&O\Big (q^{-1}|S(q,a)|\sum_{x_1\le x^{1/3}} \int_{0}^{(x-x_1^3)^{1/3}} (x-x_1^3-y^3)^{-2/3}y^2 q^{1/2+\varepsilon}\mathrm{d}y \Big).
\end{align}

Note that $S(q,a)\ll q$, then the second term of (\ref{lemma3nextnextstep}) is

\begin{equation*}
\ll \sum_{x_1\le x^{1/3}} \int_{0}^{(x-x_1^3)^{1/3}} (x-x_1^3-y^3)^{-2/3}y^2 q^{1/2+\varepsilon}\mathrm{d}y =\sum_{x_1\le x^{1/3}} (x-x_1^3)^{1/3} q^{1/2+\varepsilon} \le q^{1/2+\varepsilon}x^{2/3}.
\end{equation*}

Use integration by substitution, the integral in the first term of (\ref{lemma3nextnextstep}) is $\displaystyle\frac{1}{3}(x-x_1^3)^{2/3}B\Big(\frac{4}{3}, \frac{1}{3}\Big)$, where $B(a,b)$ represents beta function. Now repeat the argument on the sum over $x_1$, we can see that the first term is
\begin{align*}
&=\frac{1}{3}B\Big(\frac{4}{3}, \frac{1}{3}\Big)q^{-2}S(q,a)^2 \sum_{x_1\le x^{1/3}} e\Big(\frac{ax_1^3}{q}\Big)(x-x_1^3)^{2/3}\\
&=\frac{2}{3}B\Big(\frac{4}{3}, \frac{1}{3}\Big)q^{-2}S(q,a)^2 \sum_{x_1\le x^{1/3}} e\Big(\frac{ax_1^3}{q}\Big)\int_{x_1}^{x^{1/3}} (x-y^3)^{-1/3} y^2\mathrm{d}y\\
&=\frac{2}{3}B\Big(\frac{4}{3}, \frac{1}{3}\Big)q^{-2}S(q,a)^2 \int_0^{x^{1/3}} (x-y^3)^{-1/3} y^2 \Big(\sum_{x_1\le y} e\Big(\frac{ax_1^3}{q}\Big)\Big)\mathrm{d}y\\
&=\frac{2}{3}B\Big(\frac{4}{3}, \frac{1}{3}\Big)q^{-2}S(q,a)^2 \int_0^{x^{1/3}} (x-y^3)^{-1/3} y^2 (q^{-1} S(q,a)y+O(q^{1/2+\varepsilon})) \mathrm{d}y,
\end{align*}

including a main term and an error term, and the error term is

\begin{equation*}
\ll 2\int_0^{x^{1/3}} (x-y^3)^{-1/3} y^2 q^{1/2+\varepsilon} \mathrm{d}y = q^{1/2+\varepsilon}x^{2/3}.
\end{equation*}

Therefore, by (\ref{lemma3nextstep}) and (\ref{lemma3nextnextstep}), we have

\begin{equation*}
\sum_{n\le x} r_3(n)e\Big(\frac{an}{q}\Big)=\frac{2}{3}B\Big(\frac{4}{3}, \frac{1}{3}\Big)q^{-3}S(q,a)^3 \int_0^{x^{1/3}} (x-y^3)^{-1/3} y^3 \mathrm{d}y+O\Big(x^{\frac{2}{3}}q^{\frac{1}{2}+\varepsilon}\Big).
\end{equation*}

Finally, after using integration by substitution, the integral in the main term is $\displaystyle\frac{1}{3}xB\Big(\frac{4}{3},\frac{2}{3}\Big)$, and the coefficient $\displaystyle\frac{2}{3}B\Big(\frac{4}{3},\frac{1}{3}\Big)\cdot\frac{1}{3}B\Big(\frac{4}{3},\frac{2}{3}\Big)=\Gamma\Big(\frac{4}{3}\Big)^3$. The lemma follows.
\end{proof}

\vspace{3mm}

By prime factorization, every positive integer $r$ can be written uniquely as

\begin{equation}\label{fromrtor1r2r3}
r=r_1r_2^2r_3^3
\end{equation}

where $r_1$ and $r_2$ are squarefree, and $(r_1, r_2)=1$. Then the following lemma holds:

\vspace{3mm}

\textbf{Lemma 4.} \textit{If $(r,a)=1$, then for every $\varepsilon>0$,}

\begin{equation}\label{lemma4}
|S(r,a)|\ll r^{\varepsilon} r_1^{1/2} r_2 r_3^2,
\end{equation}

\textit{where $r_1$, $r_2$ and $r_3$ are given by (\ref{fromrtor1r2r3}).}

\begin{proof}
Obviously, the case that $r=1$ is correct. So it can be supposed that $r\ge2$.

\vspace{3mm}

First consider the case when $r=p^{\alpha}$, where $p$ is a prime number and $\alpha$ is a positive integer. We need to use the following result:

\vspace{3mm}

\textit{``Suppose that $p\nmid a$ and $l>\gamma$. Then $S(p^l,a)=p^{l-1}$ when $l\le k$, and $S(p^l,a)=p^{k-1}S(p^{l-k},a)$ when $l>k$."}

\vspace{3mm}

In this case $k=3$, $\gamma=2$ (if $p=3$) or $1$ (otherwise). Now let $\alpha=3u+v$, where $u$ is a nonnegative integer and $1\le v\le 3$. There are 3 different cases for $v$:

\vspace{3mm}

\underline{Case 1.} If $v=3$, then $r_1=1$, $r_2=1$, $r_3=p^{u+1}$. We have:

\begin{equation*}
S(p^{\alpha}, a)= S(p^{3u+3},a)=p^{2u+2}=r_1^{1/2}r_2r_3^2.
\end{equation*}

\vspace{3mm}

\underline{Case 2.} If $v=2$, then $r_1=1$, $r_2=p$, $r_3=p^u$. Similarly, $S(p^{\alpha},a)=S(p^{3u+2},a)=p^{2u}S(p^2,a)$.

\vspace{3mm}

If $p\ne 3$, then $S(p^2,a)=p$. So

\begin{equation*}
S(p^{\alpha},a)=p^{2u+1}=r_1^{1/2}r_2r_3^2.
\end{equation*}

If $p=3$, then an immediate calculation gives $|S(p^2,a)|\le 9$ ($=3p$). So

\begin{equation*}
|S(p^{\alpha},a)|\le 3p^{2u+1}=3r_1^{1/2}r_2r_3^2.
\end{equation*}

\vspace{3mm}

\underline{Case 3.} If $v=1$, then $r_1=p$, $r_2=1$, $r_3=p^u$. Similarly, $S(p^{\alpha},a)=S(p^{3u+1},a)=p^{2u}S(p,a)$.

\vspace{3mm}

We have $S(p,a)=\displaystyle\sum_{\chi\in\mathcal{A}} \overline{\chi}(a)\tau(\chi)$, where $\chi$ represents a character modulo $p$, $|\tau(\chi)|=p^{1/2}$ and card $\mathcal{A}=(3, p-1)-1\le 2$ (it is 2 when $p \equiv 1  (\text{mod} 3)$ and 0 when $p=3$ or $p \equiv 2  (\text{mod} 3)$). So $|S(p,a)|\le 2p^{1/2}$, therefore

\begin{equation*}
|S(p^{\alpha},a)|\le 2p^{2u+1/2}=2r_1^{1/2}r_2r_3^2.
\end{equation*}

To sum up, we always have

\begin{equation}\label{lemma4summary}
|S(p^{\alpha},a)|\le 3r_1^{1/2}r_2r_3^2.
\end{equation}

Finally, for every integer $r\ge 2$, use prime factorization: $\displaystyle r=\prod_{i=1}^tp_i^{\alpha_i}$, where $p_i$ are distinct primes and $\alpha_i$ are positive integers for $1\le i\le t$. Then by induction, there exist $a_i$ for $1\le i\le t$, such that $(p_i, a_i)=1$ and

\begin{equation*}
S(r,a)=\prod_{i=1}^t S(p_i^{\alpha_i}, a_i)
\end{equation*}

From definition, $r_1$, $r_2$ and $r_3$ are multiplicative as functions of $r$. So by (\ref{lemma4summary}),

\begin{equation*}
|S(r,a)|\le 3^t r_1^{1/2}r_2r_3^2,
\end{equation*}

where $t=\omega(r)$. Hence $3^t\ll r^{\varepsilon}$ for every $\varepsilon>0$, and the lemma follows.
\end{proof}

\vspace{3mm}

As mentioned in \S 2, the Hardy - Littlewood Method would be pretty useful if there exist $\nu(q,a)$ and $J(\beta)$, such that $\nu(q,a)J(\beta)$ is an approximation of $G(\beta+a/q)$. Fortunately, those functions exist when $a_n=r_3(n)$:

\begin{equation}\label{nuofqa}
\nu(q,a)=\frac{S(q,a)^3}{q^3},
\end{equation}

\begin{equation}\label{jofbeta}
J(\beta)=\Gamma\Big(\frac{4}{3}\Big)^3\sum_{n\le x} e(\beta n).
\end{equation}

The lemma below tells us how good the approximation is.

\vspace{3mm}

\textbf{Lemma 5.} \textit{Suppose that $(q,a)=1$ and $x\ge 1$, then for every $\varepsilon>0$,}

\begin{equation}\label{lemma5}
G\Big(\beta+\frac{a}{q}\Big)=\nu(q,a) J(\beta)+ O\Big(x^{\frac{2}{3}}q^{\frac{1}{2}+\varepsilon}(1+x|\beta|)\Big),
\end{equation}

\textit{where $G(\alpha)$, $\nu(q,a)$ and $J(\beta)$ are defined by (\ref{gofalphaforr3n}), (\ref{nuofqa}) and (\ref{jofbeta}) respectively.}

\begin{proof}
First consider the case $\beta=0$. By Lemma 3, the left-hand side of (\ref{lemma5}) is

\begin{equation*}
G\Big(\frac{a}{q}\Big)=\sum_{n\le x} r_3(n)e\Big(\frac{an}{q}\Big)=\Gamma\Big(\frac{4}{3}\Big)^3xq^{-3}S(q,a)^3+O\Big(x^{\frac{2}{3}}q^{\frac{1}{2}+\varepsilon}\Big),
\end{equation*}

while the right-hand side is

\begin{equation*}
\nu(q,a) J(0) +O\Big(x^{\frac{2}{3}}q^{\frac{1}{2}+\varepsilon}\Big)=\Gamma\Big(\frac{4}{3}\Big)^3[x]q^{-3}S(q,a)^3+O\Big(x^{\frac{2}{3}}q^{\frac{1}{2}+\varepsilon}\Big),
\end{equation*}

and the difference between the main terms is $\ll \displaystyle\Gamma\Big(\frac{4}{3}\Big)^3 q^{-3} |S(q,a)|^3 \ll 1$. Thus we obtain the lemma in the case $\beta=0$.

\vspace{3mm}

Now consider the case when $\beta$ is not necessarily 0. The left-hand side of (\ref{lemma5}) is
\begin{align*}
G\Big(\beta+\frac{a}{q}\Big)&=\sum_{n\le x} r_3(n)e\Big(\frac{an}{q}\Big)e(\beta n)\\
&=\sum_{n\le x} r_3(n) e\Big(\frac{an}{q}\Big) \Big(e(\beta x)-\int^x_n 2\pi i\beta e(\beta t)\mathrm{d}t\Big)\\
&=e(\beta x) G\Big(\frac{a}{q}\Big)-\int_0^x 2\pi i \beta e(\beta t) \Big(\sum_{n\le t} r_3(n) e\Big(\frac{an}{q}\Big)\Big)\mathrm{d}t
\end{align*}

by inverting the order of summation and integration.

\vspace{3mm}

Applying Lemma 3 again, the expression above is
\begin{align*}
&=\Gamma\Big(\frac{4}{3}\Big)^3xq^{-3}S(q,a)^3 e(\beta x)+O\Big(x^{\frac{2}{3}}q^{\frac{1}{2}+\varepsilon}\Big)\\
&-\int_0^x 2\pi i \beta e(\beta t)\Gamma\Big(\frac{4}{3}\Big)^3tq^{-3}S(q,a)^3 \mathrm{d}t +O\Big(\int_0^x |\beta| t^{\frac{2}{3}}q^{\frac{1}{2}+\varepsilon}\mathrm{d}t\Big)\\
&=\Gamma\Big(\frac{4}{3}\Big)^3 q^{-3}S(q,a)^3 \Big(xe(\beta x)-\int_0^x 2\pi i \beta e(\beta t) t \mathrm{d} t\Big)+O\Big(x^{\frac{2}{3}}q^{\frac{1}{2}+\varepsilon} (1+x|\beta|)\Big).
\end{align*}

Integrating by parts gives

\begin{equation*}
\int_0^x 2\pi i \beta e(\beta t) t \mathrm{d} t=x e(\beta x)-\int_0^x e(\beta t)\mathrm{d} t,
\end{equation*}

so
\begin{align*}
G\Big(\beta+\frac{a}{q}\Big)&=\Gamma\Big(\frac{4}{3}\Big)^3 q^{-3}S(q,a)^3 \int_0^x e(\beta t) \mathrm{d} t+O\Big(x^{\frac{2}{3}}q^{\frac{1}{2}+\varepsilon} (1+x|\beta|)\Big)\\
&=\nu(q,a) J(\beta) +\Gamma\Big(\frac{4}{3}\Big)^3\nu(q,a)\Big(\int_0^x e(\beta t) \mathrm{d} t - \sum_{n\le x} e(\beta n)\Big)+O\Big(x^{\frac{2}{3}}q^{\frac{1}{2}+\varepsilon} (1+x|\beta|)\Big).
\end{align*}

Finally, use a similar argument on $\displaystyle\int_0^x e(\beta t)\mathrm{d} t-\sum_{n\le x} e(\beta n)$, we have
\begin{align*}
\sum_{n\le x} e(\beta n)&=\sum_{n\le x} \Big(e(\beta x)-\int_n^x 2\pi i\beta e(\beta t)\mathrm{d} t\Big)\\
&=e(\beta x)\sum_{n\le x} 1-\int_0^x 2\pi i \beta e(\beta t)\Big(\sum_{n\le t} 1\Big)\mathrm{d} t\\
&=xe(\beta x) -\int_0^x 2\pi i \beta e(\beta t) t\mathrm{d} t +O(1+x|\beta|)\\
&=\int_0^x e(\beta t)\mathrm{d}t+O(1+x|\beta|)
\end{align*}

by integrating by parts. We also have $\nu(q,a)\ll 1$, so the lemma follows.
\end{proof}

\vspace{3mm}

Now define $\Delta(q,a,\beta)$ as follows:

\begin{equation}\label{deltaofqabeta}
\Delta(q,a,\beta)=G\Big(\beta+\frac{a}{q}\Big)-\nu(q,a)J(\beta).
\end{equation}

Lemma 5 says that if $(q,a)=1$ and $x\ge 1$, then for every $\varepsilon>0$,

\begin{equation}\label{deltaofqabetaproperty}
\Delta(q,a,\beta)\ll x^{2/3}q^{1/2+\varepsilon} (1+x|\beta|).
\end{equation}

\vspace{3mm}

\begin{center}

\large\textbf{5. The main term $S_3$}

\end{center}

By (\ref{s3forr3n}), Lemma 1 and (\ref{tofr}), we have

\begin{equation}
\label{FormS3}
\Gamma\Big(\frac{4}{3}\Big)^{-6}x^{-2}S_3=\sum_{q\le Q}\frac{1}{q}\sum_{r|q} rT(r)
\end{equation}
so that
\begin{align*}
\Gamma\Big(\frac{4}{3}\Big)^{-6}x^{-2}S_3&=\sum_{r\le Q}\sum_{l\le Q/r}\frac{1}{rl}\cdot rT(r)\\
&=\sum_{r\le Q} T(r)\Big(\log\frac{Q}{r}+\gamma+O\Big(\frac{r}{Q}\Big)\Big),
\end{align*}

where $\gamma$ is Euler's constant. So

\begin{equation*}
\Gamma\Big(\frac{4}{3}\Big)^{-6}x^{-2}S_3=(\log Q+\gamma)\sum_{r\le Q} T(r)-\sum_{r\le Q} T(r)\log r+O\Big(\frac{1}{Q}\sum_{r\le Q} rT(r)\Big).
\end{equation*}

By Lemma 2, $T(r)\ll r^{-2}$, so the series $\displaystyle\sum_{r=1}^{\infty} T(r)$ and $\displaystyle\sum_{r=1}^{\infty} T(r)\log r$ are convergent, and

\begin{equation*}
\sum_{r>Q} T(r) \ll \sum_{r>Q} \frac{1}{r^2} \ll \frac{1}{Q},
\end{equation*}

\begin{equation*}
\sum_{r>Q} T(r)\log r \ll \sum_{r>Q} \frac{\log r}{r^2} \ll \frac{\log Q}{Q},
\end{equation*}

\begin{equation*}
\sum_{r\le Q} rT(r) \ll \sum_{r\le Q} \frac{1}{r} \ll \log Q.
\end{equation*}

So

\begin{equation*}
\Gamma\Big(\frac{4}{3}\Big)^{-6}x^{-2}S_3=(\log Q+\gamma)\sum_{r=1}^{\infty} T(r)-\sum_{r=1}^{\infty} T(r)\log r+O\Big(\frac{\log Q}{Q}\Big).
\end{equation*}

Now define the following constants:

\begin{equation}\label{constantc0}
C_0=\Gamma\Big(\frac{4}{3}\Big)^6\sum_{r=1}^{\infty} T(r),
\end{equation}

and

\begin{equation}\label{constantc1}
C_1=\Gamma\Big(\frac{4}{3}\Big)^6\sum_{r=1}^{\infty} T(r)\log r.
\end{equation}

Therefore,

\begin{equation}\label{conclusions3}
S_3=C_0 x^2\log Q+(\gamma C_0-C_1) x^2 +O\Big(\frac{x^2 \log Q}{Q}\Big).
\end{equation}

\vspace{3mm}

\begin{center}

\large\textbf{6. The error term $S_2-S_3$}

\end{center}

We are concerned with

\begin{equation}\label{s2minuss3}
S_2-S_3=\Gamma\Big(\frac{4}{3}\Big)^3 x\sum_{q\le Q}\sum^q_{a=1} \frac{\rho(q,a)}{q^3}\Big(\sum_{\substack{n\le x \\
n\equiv a(\text{ mod } q)
}} r_3(n)-\frac{\rho(q,a)}{q^3}\Gamma\Big(\frac{4}{3}\Big)^3x \Big).
\end{equation}

By the proof of (\ref{sumrhoofqasquared}), we have:

\begin{equation*}
\sum^q_{a=1} \frac{\rho(q,a)^2}{q^6} \Gamma\Big(\frac{4}{3}\Big)^3 x=\frac{1}{q^7}\sum_{b=1}^q S(q,b)^3\overline{S(q,b)}^3 \Gamma\Big(\frac{4}{3}\Big)^3 x.
\end{equation*}

By (\ref{rhoofqa}), we have:
\begin{align*}
\sum_{a=1}^q \frac{\rho(q,a)}{q^3} \sum_{\substack{n\le x \\
n\equiv a(\text{ mod } q)
}} r_3(n)&=\sum_{a=1}^q \frac{1}{q^4}\sum_{b=1}^q e\Big(-\frac{ba}{q} \Big) S(q,b)^3 \sum_{\substack{
n\le x \\ n\equiv a(\text{ mod } q)
}} r_3(n)\\
&=\frac{1}{q^4}\sum_{b=1}^q S(q,b)^3 \sum_{a=1}^q \sum_{\substack{n\le x\\
n\equiv a(\text{ mod } q)
}} e\Big(-\frac{ba}{q}\Big) r_3(n)\\
&=\frac{1}{q^4}\sum_{b=1}^q S(q,b)^3 \sum_{n\le x} e\Big(-\frac{bn}{q}\Big) r_3(n).
\end{align*}

So
\begin{align*}
&\sum^q_{a=1} \frac{\rho(q,a)}{q^3}\Big(\sum_{\substack{n\le x \\
n\equiv a(\text{ mod } q)
}} r_3(n)-\frac{\rho(q,a)}{q^3}\Gamma\Big(\frac{4}{3}\Big)^3x \Big)\\
=&\frac{1}{q^4} \sum_{b=1}^q S(q,b)^3 \Big(\sum_{n\le x} e\Big(-\frac{bn}{q}\Big) r_3(n)-\Gamma\Big(\frac{4}{3}\Big)^3 xq^{-3}S(q,-b)^3\Big),
\end{align*}

since $\overline{S(q,b)}=S(q,-b)$.

\vspace{3mm}

By taking out the common factor $(q,b)$, we can observe that the expression above is
\begin{align*}
&=\frac{1}{q^4}\sum_{r|q}\sum^r_{\substack{a=1 \\
(a,r)=1
}}\frac{q^3}{r^3} S(r,a)^3 \Big(\sum_{n\le x} e\Big(-\frac{an}{r}\Big) r_3(n)-\Gamma\Big(\frac{4}{3}\Big)^3 xq^{-3}q^3r^{-3}S(r,-a)^3\Big)\\
&=\frac{1}{q}\sum_{r|q}\sum^r_{\substack{a=1 \\
(a,r)=1
}}\frac{S(r,a)^3}{r^3} \Big(\sum_{n\le x} e\Big(-\frac{an}{r}\Big) r_3(n)-\Gamma\Big(\frac{4}{3}\Big)^3 xr^{-3} S(r,-a)^3\Big),
\end{align*}

and note that a similar argument was used when proving Lemma 1.

\vspace{3mm}

By Lemma 3, this is
\[
\ll \frac{1}{q}\sum_{r|q}\sum^r_{\substack{a=1 \\
(a,r)=1
}}\frac{|S(r,a)|^3}{r^3} x^{2/3} r^{1/2+\varepsilon}.
\] So we have
\begin{align*}
S_2-S_3 &\ll x^{5/3}\sum_{q\le Q}\frac{1}{q}\sum_{r|q}\sum_{\substack{a=1 \\ (a,r)=1
}}^r |S(r,a)|^3 r^{-5/2+\varepsilon}\\
&=x^{5/3}\sum_{r\le Q}\Big(\sum_{l\le Q/r} \frac{1}{l}\Big) \sum_{\substack{a=1 \\ (a,r)=1
}}^r |S(r,a)|^3 r^{-7/2+\varepsilon}\\
&\ll x^{5/3} (\log Q) \sum_{r\le Q} r^{-7/2+\varepsilon} \sum^r_{\substack{a=1\\ (a,r)=1
}} |S(r,a)|^3
\end{align*}

by inserting the approximation above in (\ref{s2minuss3}).

\vspace{3mm}

By Lemma 4, if $(r,a)=1$, then for every $\varepsilon>0$, $|S(r,a)|\ll r^{\varepsilon} r_1^{1/2} r_2 r_3^2$, where $r_1$, $r_2$ and $r_3$ are defined by (\ref{fromrtor1r2r3}). So
\begin{align*}
S_2-S_3&\ll x^{5/3} (\log Q) \sum_{r\le Q} r^{-7/2+\varepsilon} \sum^r_{\substack{a=1\\ (a,r)=1
}} r^{3\varepsilon}r_1^{3/2} r_2^3 r_3^6\\
&\le x^{5/3} (\log Q) \sum_{r\le Q} r^{-5/2+4\varepsilon} r_1^{3/2} r_2^3 r_3^6\\
&= x^{5/3} (\log Q) \sum_{\substack{r_1, r_2 \text{squarefree and coprime} \\ r_1r_2^2r_3^3\le Q
}} r_1^{4\varepsilon-1} r_2^{8\varepsilon-2} r_3^{12\varepsilon-3/2}\\
&\le x^{5/3} (\log Q) \Big(\sum_{r_1\le Q} r_1^{4\varepsilon-1}\Big) \Big(\sum_{r_2=1}^{\infty} r_2^{8\varepsilon-2}\Big) \Big(\sum_{r_3=1}^{\infty} r_3^{12\varepsilon-3/2}\Big).
\end{align*}

If $\varepsilon$ is small enough, then the series above are convergent. Also,

\begin{equation*}
\sum_{r_1\le Q} r_1^{4\varepsilon-1} \le Q^{4\varepsilon}\sum_{r_1\le Q} r_1^{-1}\ll Q^{4\varepsilon}\log Q \ll Q^{5\varepsilon}.
\end{equation*}

Hence $S_2-S_3\ll x^{5/3} Q^{6\varepsilon}$. Finally, since $\varepsilon$ represents an arbitrary positive real number, then $6\varepsilon$ also has the same meaning, so we can rewrite $6\varepsilon$ as $\varepsilon$. Therefore,

\begin{equation}\label{conclusions2minuss3}
S_2-S_3 = O\Big(x^{\frac{5}{3}}Q^{\varepsilon}\Big).
\end{equation}

\vspace{3mm}

\begin{center}

\large\textbf{7. The major arcs}

\end{center}

By (\ref{deltaofqabeta}), we have
\begin{align*}
\Big|G\Big(\beta+\frac{a}{q}\Big)\Big|^2&=\Big(\nu(q,a)J(\beta)+\Delta(q,a,\beta)\Big)\Big(\overline{\nu(q,a)}\cdot\overline{J(\beta)}+\overline{\Delta(q,a,\beta)}\Big)\\
&=|\nu(q,a)|^2 |J(\beta)|^2 +\Delta_1+\Delta_2,
\end{align*}

where $\Delta_1=2\Re \Big(\overline{\nu(q,a)}\cdot\overline{J(\beta)}\cdot\Delta(q,a,\beta)\Big)$ and $\Delta_2=|\Delta(q,a,\beta)|^2$. So by (\ref{definitionofs4}), if we define

\begin{equation}\label{definitionofs5}
S_5=\sum_{q\le x/R}\int_{-1/2qR}^{1/2qR}F_q(\beta) |J(\beta)|^2 \sum_{\substack{a=1\\ (a,q)=1}}^q |\nu(q,a)|^2\mathrm{d}\beta,
\end{equation}

then

\begin{equation}\label{s4ands5}
S_4=S_5+\Sigma_1+\Sigma_2,
\end{equation}

where

\begin{equation*}
\Sigma_1=\sum_{q\le x/R}\int_{-1/2qR}^{1/2qR}F_q(\beta)\sum_{\substack{a=1\\ (a,q)=1
}}^q \Delta_1 \mathrm{d}\beta,
\end{equation*}

\begin{equation*}
\Sigma_2=\sum_{q\le x/R}\int_{-1/2qR}^{1/2qR}F_q(\beta)\sum_{\substack{a=1\\ (a,q)=1
}}^q \Delta_2 \mathrm{d}\beta.
\end{equation*}

By (\ref{estimatefsubqofalpha}) and (\ref{deltaofqabetaproperty}), we have
\begin{align*}
\Sigma_2&\ll \sum_{q\le x/R}\int_{-1/2qR}^{1/2qR}\frac{x\log(2\sqrt{x}/q)}{q+qx|\beta|} \sum_{\substack{a=1\\ (a,q)=1
}}^q \Big(x^{4/3}q^{1+2\varepsilon}(1+x|\beta|)^2\Big) \mathrm{d}\beta\\
&\ll \sum_{q\le x/R}\int_{-1/2qR}^{1/2qR} x^{7/3} (\log x) q^{1+2\varepsilon} (1+x|\beta|) \mathrm{d}\beta\\
&=x^{7/3} (\log x) \sum_{q\le x/R} q^{1+2\varepsilon}\cdot \frac{1}{2qR}\cdot\Big(2+\frac{x}{2qR}\Big)\\
&\ll x^{10/3} (\log x) R^{-2} \sum_{q\le x/R} q^{-1+2\varepsilon}\\
&\ll x^{10/3+2\varepsilon} (\log x)^2 R^{-2}.
\end{align*}

We have that if $(q,a)=1$, then $\nu(q,a)\ll q^{-1}$. For $J(\beta)$, the following properties hold: $\displaystyle J(\beta)\ll \sum_{n\le x}1\le x$, and $J(\beta)\ll \displaystyle\frac{1}{||\beta||}=\frac{1}{|\beta|}$ when $\beta\in[-1/2,1/2]\setminus\{0\}$, so $J(\beta)\ll\displaystyle\frac{x}{1+x|\beta|}$ in $[-1/2, 1/2]$. Applying the properties gives
\begin{align*}
\Sigma_1&\ll \sum_{q\le x/R}\int_{-1/2qR}^{1/2qR}\frac{x\log(2\sqrt{x}/q)}{q+qx|\beta|} \sum_{\substack{a=1\\ (a,q)=1}}^q \Big(q^{-1}\cdot\frac{x}{1+x|\beta|}\cdot x^{2/3}q^{1/2+\varepsilon}(1+x|\beta|)\Big) \mathrm{d}\beta\\
&\ll \sum_{q\le x/R} \int_{-1/2qR}^{1/2qR} x^{8/3}(\log x)q^{-1/2+\varepsilon}\cdot\frac{\mathrm{d}\beta}{1+x|\beta|}\\
&=x^{8/3} (\log x) \sum_{q\le x/R} q^{-1/2+\varepsilon}\cdot \frac{2}{x}\log\Big(1+\frac{x}{2qR}\Big)\\
&\ll x^{5/3} (\log x)^2 \sum_{q\le x/R} q^{-1/2+\varepsilon}\\
&\ll x^{13/6+\varepsilon} (\log x)^2 R^{-1/2}.
\end{align*}

Hence by (\ref{s4ands5}),

\begin{equation}\label{s4ands5property}
S_4=S_5+O\Big(x^{\frac{10}{3}+2\varepsilon}(\log x)^2 R^{-2}+x^{\frac{13}{6}+\varepsilon}(\log x)^2 R^{-\frac{1}{2}}\Big).
\end{equation}

To estimate $S_5$, note that by (\ref{nuofqa}) and (\ref{tofr}),

\begin{equation*}
\sum^q_{\substack{a=1 \\
(a,q)=1
}} |\nu(q,a)|^2=\frac{1}{q^6}\sum^q_{\substack{a=1\\
(a,q)=1
}} |S(q,a)|^6 = qT(q),
\end{equation*}

so

\begin{equation*}
S_5=\sum_{q\le x/R} qT(q) \int_{-1/2qR}^{1/2qR}F_q(\beta) |J(\beta)|^2 \mathrm{d}\beta.
\end{equation*}

Define $S_{\frac{11}{2}}$ and $S_6$ as follows:

\begin{equation*}
S_{\frac{11}{2}}=\sum_{q\le x/R} qT(q) \int_{-1/2}^{1/2}F_q(\beta) |J(\beta)|^2 \mathrm{d}\beta,
\end{equation*}

\begin{equation}\label{definitionofs6}
S_6=\sum_{q\le \sqrt{x}} qT(q) \int_{-1/2}^{1/2}F_q(\beta) |J(\beta)|^2 \mathrm{d}\beta.
\end{equation}

By (\ref{estimatefsubqofalpha}), we have a cruder estimate that $F_q(\beta)\ll q^{-1}x\log x$. By Lemma 2, we have $T(q)\ll q^{-2}$. For $J(\beta)$, we have $J(\beta)\ll \displaystyle\frac{1}{|\beta|}$ when $\beta\in[-1/2, 1/2]\setminus\{0\}$, and $\displaystyle\int_{-1/2}^{1/2} |J(\beta)|^2\mathrm{d}\beta =[x]$. So
\begin{align*}
\Big|S_{\frac{11}{2}}-S_5\Big|&\le \sum_{q\le x/R} qT(q) \Big(\int_{-1/2}^{-1/2qR}+\int_{1/2qR}^{1/2}\Big)|F_q(\beta)| |J(\beta)|^2 \mathrm{d}\beta\\
&\ll x\log x\sum_{q\le x/R} T(q) \Big(\int_{-1/2}^{-1/2qR}+\int_{1/2qR}^{1/2}\Big) \frac{1}{|\beta|^2} \mathrm{d}\beta\\
&=2x\log x\sum_{q\le x/R} T(q)(2qR-2)\\
&\ll Rx\log x \sum_{q\le x/R} q^{-1}\\
&\ll Rx(\log x)^2,
\end{align*}

and
\begin{align*}
\Big|S_6-S_{\frac{11}{2}}\Big|&\le \sum_{x/R < q\le \sqrt{x}} qT(q) \int_{-1/2}^{1/2}|F_q(\beta)| |J(\beta)|^2 \mathrm{d}\beta\\
&\ll x\log x \sum_{x/R<q\le \sqrt{x}} T(q)[x]\\
&\le x\log x\sum_{q\le \sqrt{x}} T(q)\cdot (qR)\\
&\ll Rx(\log x)^2.
\end{align*}

Hence

\begin{equation}\label{s5ands6property}
S_5=S_6+O(Rx(\log x)^2).
\end{equation}

Therefore, by (\ref{vxqwiths2s3s4forr3n}), (\ref{conclusions3}), (\ref{conclusions2minuss3}), (\ref{s4ands5property}) and (\ref{s5ands6property}),

\begin{equation}\label{vxqmaintermwiths6}
V(x,Q)=2S_6 -C_0x^2\log Q +(C_1-\gamma C_0)x^2+Q \sum_{n\le x} r_3(n)^2+U(x,Q),
\end{equation}

where

\begin{equation}\label{vxqerrortermwithR}
U(x,Q) \ll R(\log x)\sum_{n\le x}r_3(n)^2+\frac{x^{\frac{10}{3}+2\varepsilon}(\log x)^2}{R^2}+ \frac{x^{\frac{13}{6}+\varepsilon}(\log x)^2}{R^{\frac{1}{2}}}+ Rx(\log x)^2+\frac{x^2\log Q}{Q}+x^{\frac{5}{3}}Q^{\varepsilon}.
\end{equation}

\vspace{3mm}

\begin{center}

\large\textbf{8. The optimal choice for $R$}

\end{center}

First, we can redefine the upper exponent, the lower exponent and the exponent of $x$ in an expression, which are similar to the definitions given by (\ref{approximationr3nsquaredsumpositive}), (\ref{approximationr3nsquaredsumnegative}) and (\ref{approximationr3nsquaredsum}). So the choice for $R$ is optimal only if the sum of the error terms (\ref{vxqerrortermwithR}) has the minimum upper exponent of $x$, which equals to the minimum exponent of $x$ if it exists.

\vspace{3mm}

Now assume that the exponent of $x$ in $R$ exists and it equals to $B$. Since $\log x\ll x^{\varepsilon}$ for any $\varepsilon>0$, and $\varepsilon$ is not considered into the exponent, then the upper exponents of $x$ in the first four terms of (\ref{vxqerrortermwithR}) are $B+A^+$, $10/3-2B$, $13/6-B/2$ and $B+1$ respectively, where $A^+$ is defined by (\ref{approximationr3nsquaredsumpositive}) with the range given by (\ref{exponentofr3nsquaresum}). For the last two terms, we just need to consider $x^{5/3}Q^{\varepsilon}$ in which the exponent of $x$ is $5/3$, because it dominates the other one:

\begin{equation*}
\frac{x^2\log Q}{Q}\le \frac{x^2\log x}{x^{1/2}\log x}\ll x^{5/3}Q^{\varepsilon}
\end{equation*}

since $x^{1/2}\log x\le Q\le x$, which was mentioned earlier in \S 1 and \S 2. Hence the upper exponent of $x$ in (\ref{vxqerrortermwithR}) is

\begin{equation*}
I(B)=\max \Big\{B+A^+, \frac{10}{3}-2B, \frac{13}{6}-\frac{B}{2}, B+1, \frac{5}{3}\Big\}.
\end{equation*}

The function $I(B)$ has the minimum value $10/9+(2/3)A^+$ at $B=10/9-A^+/3$ when $A^+\in [1, 7/6]$. In this case, $13/18\le B\le 7/9$ and $16/9\le I(B) \le 17/9$, in other words, $R$ still satisfies (\ref{R}) if $x$ is large enough, and the sum of the error terms (\ref{vxqerrortermwithR}) is ``strictly" $\ll x^2$.

\vspace{3mm}

From the discussion above, the exponent of $x$ in the optimal choice for $R$ is $10/9-A^+/3$. But since the value of $A^+$ and the behavior of $\displaystyle\sum_{n\le x} r_3(n)^2$ are unknown, we cannot determine the expression of $R$ such that the error term is strictly optimal. Instead, we may use any expression in which the exponent of $x$ is $10/9-A^+/3$. There are two possible expressions below:

\vspace{3mm}

\underline{Possibility 1.} We may use the exponential function. In other words, let

\begin{equation}\label{exponentialR}
R=x^{\frac{10}{9}-\frac{A^+}{3}}.
\end{equation}

By (\ref{approximationr3nsquaredsum}), $\displaystyle\sum_{n\le x} r_3(n)^2\ll x^{(A^+)+\varepsilon}$ for any $\varepsilon>0$, and also note that $\log x\ll x^{\varepsilon}$, then all the other terms in (\ref{vxqerrortermwithR}) are $\ll$ the second term, which is $\ll x^{10/9+(2/3)(A^+)+4\varepsilon}$. Since $4\varepsilon$ represents an arbitrary positive real number and $\varepsilon$ has the same meaning, we can rewrite $4\varepsilon$ as $\varepsilon$. Therefore, the error term is

\begin{equation}\label{errorexp}
U(x,Q)\ll x^{\frac{10}{9}+\frac{2}{3}(A^+)+\varepsilon}.
\end{equation}

\vspace{3mm}

\underline{Possibility 2.} If the exponent of $x$ in $\displaystyle\sum_{n\le x} r_3(n)^2$ exists, then $A=A^+$, and the exponent of $x$ in $R$ is $10/9-A/3$. In this case, we can avoid using $A$ or $A^+$, and instead, we may use $\displaystyle\sum_{n\le x} r_3(n)^2$. Let

\begin{equation}\label{r3nR}
R=x^{\frac{10}{9}}\Big(\sum_{n\le x} r_3(n)^2\Big)^{-\frac{1}{3}},
\end{equation}

then again, all the other terms in (\ref{vxqerrortermwithR}) are $\ll$ the second term since $\displaystyle x\ll \sum_{n\le x}r_3(n)^2$. Use a similar argument, we can see that the error term is

\begin{equation}\label{errorr3n}
U(x,Q)\ll x^{\frac{10}{9}+\varepsilon} \Big(\sum_{n\le x} r_3(n)^2\Big)^{\frac{2}{3}}.
\end{equation}

Right now, either (\ref{errorexp}) or (\ref{errorr3n}) is a possible choice for the approximation, so both of them will be included in the conclusion. If the behavior of $\displaystyle\sum_{n\le x} r_3(n)^2$ is found in the future, then we may find the $R$ which is strictly optimal. However, at least now we have the result regarding the exponent of $x$ in the error term, in other words, we know how ``large" the error term is.

\vspace{3mm}

\begin{center}

\large\textbf{9. The main term $S_6$}

\end{center}

We start from the integral in (\ref{definitionofs6}). By (\ref{fsubqofalpha}) and (\ref{jofbeta}),

\begin{equation*}
\int_{-1/2}^{1/2} F_q(\beta)|J(\beta)|^2 \mathrm{d}\beta=\Gamma\Big(\frac{4}{3}\Big)^6\sum_{\substack{l\le\sqrt{x}\\ q|l
}}\Big(\sum_{m\le x/l}+\sum_{\sqrt{x}<m\le \min (Q, x/l)}\Big) ([x]-lm),
\end{equation*}

and a straightforward calculation shows that

\begin{equation*}
\Big(\sum_{m\le x/l}+\sum_{\sqrt{x}<m\le \min (Q, x/l)}\Big) ([x]-lm)
\end{equation*}

is

\begin{equation*}
\frac{x^2}{2l}+\frac{x}{2l}(\sqrt{x}-l)^2-\frac{Q^2}{2l}\Big(\frac{x}{Q}-l\Big)^2+O(x)
\end{equation*}

when $l\le x/Q$, and is

\begin{equation*}
\frac{x^2}{2l}+\frac{x}{2l}(\sqrt{x}-l)^2+O(x)
\end{equation*}

when $x/Q<l\le \sqrt{x}$. For convenience, we write it as $K(x,l,Q)$, then by (\ref{definitionofs6}),

\begin{equation*}
S_6=\Gamma\Big(\frac{4}{3}\Big)^6\sum_{q\le \sqrt{x}}qT(q)\sum_{\substack{l\le\sqrt{x} \\
q|l
}} K(x,l,Q)=\Gamma\Big(\frac{4}{3}\Big)^6\sum_{l\le \sqrt{x}}\Big(\sum_{q|l} qT(q)\Big) K(x,l,Q).
\end{equation*}

Define

\begin{equation}\label{hofl}
h(l)=\sum_{q|l} qT(q).
\end{equation}

Then
\begin{align*}
S_6&=\Gamma\Big(\frac{4}{3}\Big)^6\Big(\sum_{l\le x/Q}h(l) K(x,l,Q)+ \sum_{x/Q<l\le \sqrt{x}} h(l)K(x,l,Q)\Big)\\
&=\Gamma\Big(\frac{4}{3}\Big)^6\Big(\frac{x^2}{2}\sum_{l\le\sqrt{x}}\frac{h(l)}{l}+\frac{x}{2}\sum_{l\le \sqrt{x}}\frac{h(l)}{l}(\sqrt{x}-l)^2-\frac{Q^2}{2}\sum_{l\le x/Q}\frac{h(l)}{l}\Big(\frac{x}{Q}-l\Big)^2+O\Big(x\sum_{l\le\sqrt{x}}h(l)\Big).
\end{align*}

Actually, we are more interested in $2S_6$ than $S_6$, so if we define

\begin{equation}\label{wofX}
W(X)=\sum_{l\le X} \frac{h(l)}{l}(X-l)^2,
\end{equation}

then

\begin{equation}\label{2s6initial}
2S_6=\Gamma\Big(\frac{4}{3}\Big)^6\Big(x^2\sum_{l\le\sqrt{x}}\frac{h(l)}{l}+xW(\sqrt{x})-Q^2W\Big(\frac{x}{Q}\Big)\Big)+O\Big(x\sum_{l\le\sqrt{x}}h(l)\Big).
\end{equation}

By Lemma 2 and (\ref{hofl}), the error term is

\begin{equation*}
\ll x\sum_{l\le\sqrt{x}}\sum_{q|l}qT(q)=x\sum_{q\le \sqrt{x}}qT(q)\sum_{\substack{l\le\sqrt{x}\\
q|l
}} 1\le x^{3/2}\sum_{q\le\sqrt{x}}T(q)\ll x^{3/2}.
\end{equation*}

The first main term is

\begin{equation*}
\Gamma\Big(\frac{4}{3}\Big)^6 x^2\sum_{l\le\sqrt{x}}\frac{h(l)}{l}=\Gamma\Big(\frac{4}{3}\Big)^6 x^2\sum_{l\le\sqrt{x}}\frac{1}{l}\sum_{q|l} qT(q),
\end{equation*}

which looks similar to the formula for $S_3$ in (\ref{FormS3})

\begin{equation*}
S_3=\Gamma\Big(\frac{4}{3}\Big)^6 x^2\sum_{q\le Q}\frac{1}{q} \sum_{r|q} rT(r)
\end{equation*}

that was mentioned in \S 5. So we can substitute $Q$, $q$ and $r$ by $\sqrt{x}$, $l$ and $q$ respectively from the conclusion (\ref{conclusions3}), and we have:
\begin{align}\label{2s6firstmainterm}
\Gamma\Big(\frac{4}{3}\Big)^6 x^2\sum_{l\le\sqrt{x}} \frac{h(l)}{l}&=C_0 x^2\log\sqrt{x} +(\gamma C_0-C_1)x^2+O\Big(\frac{x^2\log\sqrt{x}}{\sqrt{x}}\Big)\notag\\
&=\frac{1}{2}C_0 x^2\log x +(\gamma C_0-C_1)x^2+O(x^{3/2}\log\ x),
\end{align}

where the constants $C_0$ and $C_1$ are defined by (\ref{constantc0}) and (\ref{constantc1}). So by (\ref{2s6initial}) and (\ref{2s6firstmainterm}),

\begin{equation}\label{2s6withwofXonly}
2S_6=\frac{1}{2}C_0 x^2\log x +(\gamma C_0-C_1)x^2+\Gamma\Big(\frac{4}{3}\Big)^6\Big(xW(\sqrt{x})-Q^2W\Big(\frac{x}{Q}\Big)\Big)+O(x^{3/2}\log x).
\end{equation}

So in order to draw to a conclusion, we only need to evaluate $W(\sqrt{x})$ and $W\Big(\displaystyle\frac{x}{Q}\Big)$. From now on we will focus on evaluating $W(X)$, where $1\le X\le \sqrt{x}$.

\vspace{3mm}

\begin{center}

\large\textbf{10. The final result}

\end{center}

To evaluate $W(X)$, we need the following lemmas:

\vspace{3mm}

\textbf{Lemma 6.} \textit{$T(r)$ is a multiplicative function}.

\begin{proof}
Recall that $T(r)$ is defined by (\ref{tofr}). Obviously $T(1)=1$, so it is sufficient to prove that if $(r_1, r_2)=1$, then $T(r_1r_2)=T(r_1)T(r_2)$.

\vspace{3mm}

By definition,

\begin{equation*}
T(r_1r_2)=\frac{1}{r_1^7 r_2^7}\sum_{\substack{c=1\\
(c, r_1 r_2)=1
}}^{r_1 r_2} |S(r_1, cr_2^2)|^6 |S(r_2, cr_1^2)|^6.
\end{equation*}

From elementary number theory, we know that if $1\le c\le r_1r_2$, $(c, r_1r_2)=1$ and $(r_1, r_2)=1$, then there exists a unique pair of numbers $c_1$ and $c_2$, such that $c \equiv c_2r_1+c_1r_2 (\text{mod} r_1r_2)$, $1\le c_i\le r_i$ and $(c_i, r_i)=1$ where $i=1$ or $2$. In this situation, we have $S(r_1, cr_2^2)=S(r_1, c_1)$ and $S(r_2, cr_1^2)=S(r_2, c_2)$ by definition (\ref{sofqb}). So the identity above becomes

\begin{equation*}
T(r_1 r_2)=\frac{1}{r_1^7 r_2^7}\sum_{\substack{c_1=1 \\
(c_1, r_1)=1
}}^{r_1}\sum_{\substack{c_2=1\\
(c_2, r_2)=1
}}^{r_2} |S(r_1, c_1)|^6 |S(r_2, c_2)|^6 =T(r_1)T(r_2),
\end{equation*}

which completes the proof.
\end{proof}

\vspace{3mm}


\textbf{Lemma 7.} \textit{Assume that $s=\sigma+it$ is a complex number, where $\sigma, t\in \mathbb{R}$ and $\sigma>-2$. Then the following results hold:}

\vspace{3mm}

\textit{(1). If $p\ne 3$ is a prime number, then}

\begin{equation*}
\sum_{k=0}^{\infty} \frac{T(p^k)}{p^{ks}} = \frac{1}{1-p^{-(3s+6)}} \Big(1+\frac{1}{p^{s+7}}\sum_{c=1}^{p-1} |S(p,c)|^6 +\frac{p-1}{p^{2s+7}}-\frac{1}{p^{3s+7}}\Big).
\end{equation*}

\textit{(2).}

\begin{equation*}
\sum_{k=0}^{\infty} \frac{T(3^k)}{3^{ks}} = \frac{1}{1-3^{-(3s+6)}} \Big(1+\frac{1}{3^{2s+14}}\sum_{\substack{c=1 \\
(c,3)=1
}}^{9} |S(9,c)|^6 -\frac{1}{3^{3s+7}}\Big).
\end{equation*}

\begin{proof}
By Lemma 4, if $p$ is a prime number, $(p,a)=1$ and $u$ is a nonnegative integer, then

\vspace{3mm}

(1). $S(p^{3u+3}, a)= p^{2u+2}$,

\vspace{3mm}

(2). $S(p^{3u+2}, a)=p^{2u+1}$ when $p\ne3$ and $S(3^{3u+2}, a)=3^{2u}S(9,a)$,

\vspace{3mm}

(3). $S(p^{3u+1}, a)=p^{2u} S(p,a)$.

\vspace{3mm}

In particular $S(3^{3u+1},a)=0$ since $S(3,a)=0$.

\vspace{3mm}

So by definition (\ref{tofr}), we have

\vspace{3mm}

(1). $\displaystyle T(p^{3u+3})=\frac{p-1}{p^{6u+7}}$,

\vspace{3mm}

(2). $\displaystyle T(p^{3u+2})=\frac{p-1}{p^{6u+7}}$ when $p\ne 3$ and
\[
T(3^{3u+2})=\frac{1}{3^{6u+14}}\sum_{\substack{c=1 \\
(c,3)=1
}}^{9} |S(9,c)|^6,
\]

\vspace{3mm}

(3). $\displaystyle T(p^{3u+1})=\frac{1}{p^{6u+7}}\sum_{c=1}^{p-1} |S(p,c)|^6$, therefore $T(3^{3u+1})=0$ since $S(3,1)=S(3,2)=0$.

\vspace{3mm}

Also, the identity

\begin{equation}\label{lemma7series}
\sum_{k=0}^{\infty}\frac{T(p^{k})}{p^{ks}}=1+\sum_{u=0}^{\infty}\frac{T(p^{3u+1})}{p^{(3u+1)s}}+\sum_{u=0}^{\infty}\frac{T(p^{3u+2})}{p^{(3u+2)s}}+\sum_{u=0}^{\infty}\frac{T(p^{3u+3})}{p^{(3u+3)s}}
\end{equation}

holds when the series on the right-hand side are convergent. First assume that $p\ne 3$, then from the discussion above, the right-hand side of (\ref{lemma7series}) becomes

\begin{equation*}
1+\Big(\frac{1}{p^{s+7}}\sum_{c=1}^{p-1} |S(p,c)|^6+\frac{p-1}{p^{2s+7}}+\frac{p-1}{p^{3s+7}}\Big)\sum_{u=0}^{\infty} \frac{1}{p^{(3s+6)u}}.
\end{equation*}

If $\sigma=\Re s>-2$, then the series above is absolutely convergent since $|p^{3s+6}|=p^{3\sigma+6}>1$, and it converges to $\displaystyle\frac{1}{1-p^{-(3s+6)}}$. Hence the first result is proved by (\ref{lemma7series}).

\vspace{3mm}

Now assume that $p=3$, then the right-hand side of (\ref{lemma7series}) becomes

\begin{equation*}
1+\Big(\frac{1}{3^{2s+14}}\sum_{\substack{c=1\\
(c,3)=1
}}^9 |S(9,c)|^6+\frac{2}{3^{3s+7}}\Big)\sum_{u=0}^{\infty} \frac{1}{3^{(3s+6)u}}.
\end{equation*}

Similarly, the series above converges to $\displaystyle\frac{1}{1-3^{-(3s+6)}}$ when $\sigma>-2$, thus the second result is proved.
\end{proof}

\vspace{3mm}

Now let

\begin{equation}\label{Dofs}
D(s)=\sum_{l=1}^{\infty} \frac{h(l)}{l}\cdot\frac{1}{l^s}=\sum_{l=1}^{\infty}\frac{h(l)}{l^{s+1}}
\end{equation}

if $\sigma=\Re s>0$. Then by (\ref{hofl}), Lemma 6 and Lemma 7, we have
\begin{align}
D(s)&=\sum_{l=1}^{\infty} \frac{1}{l^{s+1}} \sum_{q|l} qT(q)\notag\\
&=\sum_{q=1}^{\infty}\sum_{m=1}^{\infty}\frac{qT(q)}{q^{s+1} m^{s+1}}\notag\\
&=\zeta(s+1)\sum_{q=1}^{\infty} \frac{T(q)}{q^s} \label{Dofsstep1}\\
&=\zeta(s+1)\prod_p \Big(\sum_{k=0}^{\infty} \frac{T(p^k)}{p^{ks}} \Big)\notag\\
&=\zeta(s+1) \zeta(3s+6) D_0(s) \label{Dofsstep2},
\end{align}

where $D_0(s)$ is defined as

\begin{equation}\label{D0ofs}
\prod_{p\ne 3} \Big(1+\frac{1}{p^{s+7}}\sum_{c=1}^{p-1} |S(p,c)|^6 +\frac{p-1}{p^{2s+7}}-\frac{1}{p^{3s+7}}\Big)\cdot \Big(1+\frac{1}{3^{2s+14}}\sum_{\substack{c=1 \\ (c,3)=1
}}^{9} |S(9,c)|^6 -\frac{1}{3^{3s+7}}\Big).
\end{equation}

By Lemma 4, we have (I) $|S(p,c)|\le 2p^{1/2}$ when $p$ is a prime and $1\le c\le p-1$, and (II) $|S(9,c)|\le 9$ when $c=1, 2, 4, 5, 7$ or $8$. So

\begin{equation*}
|D_0(s)|\le \prod_{p\ne 3} \Big(1+\frac{64}{p^{\sigma+3}} +\frac{1}{p^{2\sigma+6}}+\frac{1}{p^{3\sigma+7}}\Big)\cdot \Big(1+\frac{2}{3^{2\sigma+1}} +\frac{1}{3^{3\sigma+7}}\Big).
\end{equation*}

The product on the right-hand side of the inequality above converges locally uniformly when $\sigma>-2$, so $D_0(s)$ converges absolutely and locally uniformly for $\Re s>-2$. Note that Lemma 7 also holds when $\Re s>-2$, so (\ref{Dofsstep2}) affords an analytic continuation of $D(s)$ to that open half plane. Thus from the properties of the zeta function, $D(s)$ is meromorphic in $\Re s>-2$ with simple poles at $s=0$ (from $s+1=1$) and $s=-5/3$ (from $3s+6=1$).

\vspace{3mm}

Moreover, we have the following identity if we compare (\ref{Dofsstep1}) with (\ref{Dofsstep2}) without considering the analytic continuation:

\begin{equation}\label{DofsCor}
\sum_{q=1}^{\infty} \frac{T(q)}{q^s}= \zeta(3s+6)D_0(s).
\end{equation}

The identity is correct when the series converges, especially when $\Re(3s+6)>1$ and $\Re s>-2$. For example, $s=-1$ satisfies both conditions, therefore:

\vspace{3mm}


\textbf{Lemma 8.} \textit{The series $\displaystyle\sum_{q=1}^{\infty} qT(q)$ converges to $\zeta(3) D_0(-1)$}.

\vspace{3mm}

Note that we cannot prove Lemma 8 by Lemma 2, since the harmonic series is divergent.

\vspace{3mm}

Now consider the relationship between $W(X)$ and $D(s)$.  W(X) is a essentiall a Ces\`{a}ro mean of the arithmetical function $h(l)/l$.  Thus as in Montgomery and Vaughan \cite{MV}, page 142 we have

\begin{equation}\label{wofxandDofs}
\frac{1}{2} W(X)= \frac{1}{2\pi i}\int_{\sigma_0-i\infty}^{\sigma_0+i\infty} D(s)\cdot\frac{X^{s+2}}{s(s+1)(s+2)}\mathrm{d}s,
\end{equation}

where $\sigma_0>0$. The integrand above is homomorphic in the open half plane $\Re s>-2$ except for simple poles at both $s=-1$ and $s=-5/3$, and a double pole at $s=0$.

\vspace{3mm}

To estimate the integral, let $0<\delta<1/3$, then $-2<-2+\delta<-5/3$. Then let $\mathcal{C}(t_1, t_2)$ denote the rectangular contour with vertices $\sigma_0-it_1$, $\sigma_0+it_2$, $-2+\delta+it_2$ and $-2+\delta-it_1$ where $t_1, t_2>1$. By the residue formula, we have:

\begin{equation}\label{residueformula}
\frac{1}{2\pi i}\int_{\mathcal{C}(t_1, t_2)} D(s)\cdot\frac{X^{s+2}}{s(s+1)(s+2)}\mathrm{d}s=\text{Res}(f,0)+\text{Res}(f, -1)+\text{Res}\Big(f, -\frac{5}{3}\Big),
\end{equation}

where $f$ represents the integrand. The residues are given by

\begin{equation*}
\text{Res}(f, 0)= \lim_{s\to 0} \frac{\mathrm{d}}{\mathrm{d}s} \Big(\frac{s^2D(s)X^{s+2}}{s(s+1)(s+2)}\Big) = D_1X^2\log X+D_2 X^2,
\end{equation*}

\begin{equation*}
\text{Res}(f, -1)= \lim_{s\to -1} \frac{(s+1)D(s)X^{s+2}}{s(s+1)(s+2)}=D_3X,
\end{equation*}

\begin{equation*}
\text{Res}\Big(f, -\frac{5}{3}\Big)= \lim_{s\to -5/3} \frac{(s+5/3)D(s)X^{s+2}}{s(s+1)(s+2)}=D_4X^{1/3},
\end{equation*}

where $D_1, D_2, D_3$ and $D_4$ are constants. By (\ref{Dofsstep2}), (\ref{DofsCor}) (let $s=0$) and Lemma 8,

\begin{equation*}
D_1= \lim_{s\to 0} \frac{sD(s)}{(s+1)(s+2)} = \frac{1}{2}\zeta(6)D_0(0) \lim_{s\to0} s\zeta(s+1) = \frac{1}{2}\sum_{q=1}^{\infty} T(q),
\end{equation*}

\begin{equation*}
D_3= \lim_{s\to -1} \frac{D(s)}{s(s+2)} = -D(-1) = -\zeta(0)\zeta(3)D_0(-1)= \frac{1}{2}\sum_{q=1}^{\infty} qT(q),
\end{equation*}

\begin{equation*}
D_4=\lim_{s\to -5/3} \frac{(s+5/3)D(s)}{s(s+1)(s+2)}=  \frac{9}{10}\zeta\Big(-\frac{2}{3}\Big) D_0\Big(-\frac{5}{3}\Big),
\end{equation*}

so the right-hand side of (\ref{residueformula}) is

\begin{equation}\label{residueformulaRHS}
\frac{1}{2}\Big(\sum_{q=1}^{\infty}T(q)\Big) X^2\log X + D_2X^2 + \frac{1}{2}\Big(\sum_{q=1}^{\infty}qT(q)\Big) X + \frac{9}{10}\zeta\Big(-\frac{2}{3}\Big) D_0\Big(-\frac{5}{3}\Big) X^{1/3}.
\end{equation}

Note that we don't need to evaluate $D_2$ since the related terms will be cancelled in the end.

\vspace{3mm}

Now we focus on the left-hand side of (\ref{residueformula}). Rewrite the integrand as $\displaystyle\frac{\zeta(s+1)\zeta(3s+6)D_0(s)X^{s+2}}{s(s+1)(s+2)}$ by (\ref{Dofsstep2}), where $D_0(s)$ is uniformly bounded when $\Re s\ge -2+\delta$. By Section 5.1 of Titchmarsh \cite{ET}, if $w$ is a complex number, then for any $\varepsilon>0$,

\begin{equation*}
\zeta(w)-\frac{1}{w-1}\ll \tau^{\lambda(u)+\varepsilon},
\end{equation*}

where $u=\Re w$, $v=\Im w$, $\tau=4+|v|$, and $\lambda(u)$ is defined as 0 when $u>1$, $1/2-(1/2)u$ when $0<u\le 1$, and $1/2-u$ when $u\le0$.

\vspace{3mm}

Consider the integral along the line segment $\sigma_0+it_2 \to -2+\delta +it_2$. If $s$ lies on the line segment, then $\Re(s+1)\ge -1+\delta>-1$ and $\Re(3s+6)\ge 3\delta>0$. Also, both $|(s+1)-1|^{-1}$ and $|(3s+6)-1|^{-1}$ are less than 1 since $\Im s=t_2>1$. So $\zeta(s+1)\ll (4+t_2)^{\lambda(\Re s+1)+\varepsilon}\ll (4+t_2)^{3/2+\varepsilon}$, and $\zeta(3s+6)\ll (4+3t_2)^{\lambda(3\Re s+6)+\varepsilon} \ll (4+3t_2)^{1/2+\varepsilon}$. So the integrand is $\ll (4+t_2)^{2+2\varepsilon}t_2^{-3} X^{\sigma_0+2}$. Hence the integral along this line segment is $\to 0$ when $t_2\to \infty$ and $\varepsilon$ is small enough, since the line segment is of finite length. Similarly, the integral along the line segment $-2+\delta-it_1 \to \sigma_0-it_1$ is $\to 0$ when $t_1\to\infty$.

\vspace{3mm}

Now consider the integral along the line segment $-2+\delta+it_2 \to -2+\delta-it_1$. By substitution, the integral (with coefficient) is

\begin{equation*}
-\frac{1}{2\pi} X^{\delta} \int_{-t_1}^{t_2} \frac{\zeta(-1+\delta+it)\zeta(3\delta+3it)D_0(-2+\delta+it)X^{it}}{(-2+\delta+it)(-1+\delta+it)(\delta+it)}\mathrm{d}t = -\frac{1}{2\pi} X^{\delta} (I_1+I_2+I_3),
\end{equation*}

where $I_1$, $I_2$ and $I_3$ are integrals with the same integrand as left-hand side but the intervals of integration are $[-t_1, -1]$, $[-1,1]$ and $[1, t_2]$ respectively. Note that $-1<-1+\delta<-2/3$ and $0<3\delta<1$ when $0<\delta<1/3$, so $\zeta(-1+\delta+it)\ll (4+|t|)^{\lambda(-1+\delta)+\varepsilon} = (4+|t|)^{3/2-\delta+\varepsilon}$ and $\zeta(3\delta+3it)\ll (4+3|t|)^{\lambda(3\delta)+\varepsilon} = (4+3|t|)^{1/2-3\delta/2+\varepsilon}$ if $t=\Im s$ satisfies $|t|\ge 1$. Also, $|X^{it}|=1$, and $D_0(-2+\delta+it)$ is bounded. So the integrand above is $\ll |t|^{3/2-\delta+\varepsilon+1/2-3\delta/2+\varepsilon-3} = |t|^{-1-5\delta/2+2\varepsilon}$. We can choose $\varepsilon=\delta>0$, then the integrand is $\ll |t|^{-1-\delta/2}$ when $|t|\ge 1$. So for a fixed $\delta$, $I_1\ll 1-t_1^{-\delta/2}$ and $I_3\ll 1-t_2^{-\delta/2}$. For $I_2$, we have $I_2\ll 1$ since the integrand above is bounded when $|t|\le 1$. Hence the integral along the line segment $-2+\delta+it_2 \to -2+\delta-it_1$ is $\to O(X^{\delta})$ when $t_1\to\infty$ and $t_2\to\infty$.

\vspace{3mm}

Finally, by (\ref{wofxandDofs}), the integral along the line segment $\sigma_0-it_1\to \sigma_0+it_2$ with coefficient is $\to W(X)/2$ when $t_1\to \infty$ and $t_2\to \infty$. So let $t_1, t_2\to \infty$ on both sides of (\ref{residueformula}), we have

\begin{equation}\label{WofXanswer}
W(X)=\Big(\sum_{q=1}^{\infty}T(q)\Big) X^2\log X + 2D_2X^2 + \Big(\sum_{q=1}^{\infty}qT(q)\Big) X + \frac{9}{5}\zeta\Big(-\frac{2}{3}\Big) D_0\Big(-\frac{5}{3}\Big) X^{1/3} +O(X^{\delta})
\end{equation}

by (\ref{residueformulaRHS}). Now substitute $X=\sqrt{x}$ into (\ref{WofXanswer}), and note that $x^{3/2}$, $x^{7/6}$ and $O(x^{1+\delta/2})$ are all $\ll x^{3/2}$, we have

\begin{equation}\label{xWofsqrtofx}
xW(\sqrt{x})= \frac{1}{2}\Big(\sum_{q=1}^{\infty}T(q)\Big) x^2\log x + 2D_2 x^2 + O(x^{3/2}).
\end{equation}

Similarly, we can get the expression of $Q^2W(x/Q)$ by substituting $X=x/Q$ into (\ref{WofXanswer}). Then by (\ref{2s6withwofXonly}) and (\ref{constantc0}),

\begin{equation}\label{2s6answer}
2S_6=C_0 x^2\log Q +(\gamma C_0-C_1) x^2 - A_1 Qx + A_2 Q^{5/3} x^{1/3} + O(x^{3/2} \log x) +O(Q^{2-\delta} x^{\delta}),
\end{equation}

where

\begin{equation}\label{A1}
A_1=\Gamma\Big(\frac{4}{3}\Big)^6 \sum_{q=1}^{\infty} qT(q),
\end{equation}

\begin{equation}\label{A2}
A_2=-\frac{9}{5} \Gamma\Big(\frac{4}{3}\Big)^6 \zeta\Big(-\frac{2}{3}\Big) D_0\Big(-\frac{5}{3}\Big).
\end{equation}

Note that both $A_1$ and $A_2$ are positive.

\vspace{3mm}

Finally, by (\ref{vxqmaintermwiths6}), (\ref{errorexp}) and (\ref{errorr3n}), and some constant-related definitions (\ref{sofqb}), (\ref{tofr}) and (\ref{D0ofs}), the final result can be stated as follows:

\vspace{3mm}


\textbf{Theorem 1.} \textit{Assume that $\varepsilon$ and $\delta$ are positive numbers satisfying $0<\delta<\displaystyle\frac{1}{3}$ and let}

\begin{equation*}
U_0(x,Q)=V(x,Q)-Q\sum_{n\le x} r_3(n)^2 +A_1 Qx -A_2 Q^{\frac{5}{3}} x^{\frac{1}{3}}
\end{equation*}

\textit{where}

\begin{equation*}
A_1=\Gamma\Big(\frac{4}{3}\Big)^6 \sum_{q=1}^{\infty} \frac{1}{q^6} \sum^q_{\substack{c=1 \\
(c,q)=1
}} |S(q,c)|^6,
\end{equation*}

\begin{equation*}
A_2=-\frac{9}{5} \Gamma\Big(\frac{4}{3}\Big)^6 \zeta\Big(-\frac{2}{3}\Big) D_0\Big(-\frac{5}{3}\Big),
\end{equation*}

\begin{equation*}
D_0(s)= \prod_{p\ne 3} \Big(1+\frac{1}{p^{s+7}}\sum_{c=1}^{p-1} |S(p,c)|^6 +\frac{p-1}{p^{2s+7}}-\frac{1}{p^{3s+7}}\Big)\cdot \Big(1+\frac{1}{3^{2s+14}}\sum_{\substack{c=1 \\
(c,3)=1
}}^{9} |S(9,c)|^6 -\frac{1}{3^{3s+7}}\Big),
\end{equation*}

\textit{and}

\begin{equation*}
S(q,c)= \sum_{m=1}^q e\Big(\frac{cm^3}{q}\Big).
\end{equation*}

\textit{Then when $x^{\frac{1}{2}}\log x \le Q \le x$ one has}

\begin{equation*}
U_0(x,Q) \ll x^{\frac{10}{9}+\frac{2}{3}(A^+)+\varepsilon}+Q^{2-\delta}x^{\delta}
\end{equation*}

\textit{where $A^+$ is the upper exponent of $x$ in the representation of $\displaystyle\sum_{n\le x} r_3(n)^2$, or}

\begin{equation*}
U_0(x,Q) \ll x^{\frac{10}{9}+\varepsilon} \Big(\sum_{n\le x} r_3(n)^2\Big)^{\frac{2}{3}}+Q^{2-\delta}x^{\delta}
\end{equation*}

\textit{if the exponent of $x$ in the representation of $\displaystyle\sum_{n\le x} r_3(n)^2$ exists.}

\vspace{3mm}

\begin{center}

\large\textbf{11. Special cases}

\end{center}

Theorem 1 gives us a general result of the problem. However, if $Q$ is very close to $x$, say $Q\asymp x$, then $Qx\asymp Q^{5/3}x^{1/3} \asymp Q^{2-\delta}x^{\delta} \asymp x^2$, so the exponent of $x$ in the error term $U_0(x,Q)$ is the same as that in part of the main terms $-A_1Qx+A_2Q^{5/3}x^{1/3}$, which is 2. So in this case, we need to find more precise results.

\vspace{3mm}

Note that the error term $O(Q^{2-\delta}x^{\delta})$ comes from $Q^2 W(x/Q)$, so instead of applying (\ref{WofXanswer}) to approximate this term we evaluate it exactly. If $k\ge 1$ is an integer, then by (\ref{wofX}), when $x/(k+1) <Q\le x/k$, we have

\begin{equation}\label{WofxoverQspecial}
Q^2 W\Big(\frac{x}{Q}\Big)=Q^2\sum_{l\le x/Q} \frac{h(l)}{l} \Big(\frac{x}{Q}-l\Big)^2=\sum_{l\le k} \frac{h(l)}{l} (x-lQ)^2,
\end{equation}

which is a sum of $k$ terms. We will apply this result to a theorem and apply the following special cases to several corollaries:

\vspace{3mm}

(I). By (\ref{hofl}), (\ref{tofr}) and (\ref{sofqb}), we have $h(1)=1$ and $h(2)=1$. So $Q^2W(x/Q)=x^2-2Qx+Q^2$ when $x/2<Q\le x$, and it is $3x^2/2-4Qx+3Q^2$ when $x/3<Q\le x/2$.

\vspace{3mm}

(II). From (I), if $Q=x$ (the largest possible value of $Q$), then $Q^2W(x/Q)=x^2W(1)=0$.

\vspace{3mm}

(III). If $Q=x/m$ where $m\ge 1$ is a constant, then $k=[m]$. So by (\ref{WofxoverQspecial}),

\begin{equation*}
Q^2 W\Big(\frac{x}{Q}\Big)=\frac{x^2}{m^2} W(m)=\sum_{l\le [m]} \frac{h(l)}{l} \Big(x-\frac{lx}{m}\Big)^2=\frac{x^2}{m^2}\sum_{l\le m} \frac{h(l)}{l} (m-l)^2.
\end{equation*}

For $xW(\sqrt{x})$, we may still apply (\ref{xWofsqrtofx}). However, in this case the term $2D_2x^2$ will not be cancelled, so we need to evaluate $D_2$, which is not easy if we evaluate the corresponding residue. So instead of using the method in \S 10, we would rather use a cruder estimate for $W(X)$. First, we state the following result as a lemma, which is proved in Vaughan \cite{Va4}.

\vspace{3mm}


\textbf{Lemma 9.} \textit{If $Y\ge 1$, then}

\begin{equation}\label{lemma9}
\sum_{m\le Y} \frac{1}{m} (Y-m)^2= Y^2\log Y+C_2Y^2+Y+O(1),
\end{equation}

where

\begin{equation}\label{constantc2}
C_2=-\frac{11}{12}-2\int_1^{\infty} \frac{B_2(u)}{u^3}\mathrm{d}u
\end{equation}

\textit{and}

\begin{equation}\label{B2ofu}
B_2(u)=\frac{1}{2} (u-[u])^2 -\frac{1}{2}(u-[u])+\frac{1}{12}.
\end{equation}

\vspace{3mm}

In order to prove the lemma, we start from proving that

\begin{equation*}
\int^Y_1 \Big(\frac{Y^2}{u^2}-1\Big) [u]\mathrm{d}u= \sum_{m\le Y} \frac{1}{m} (Y-m)^2
\end{equation*}

by evaluating the integral exactly. Then we integrate it by parts and use the fact that $B_2(u)$ is bounded to get the conclusion.

\vspace{3mm}

Now we can evaluate $W(\sqrt{x})$. By (\ref{wofX}), (\ref{hofl}), Lemma 9 (when $Y=\sqrt{x}/q$) and Lemma 2, we have
\begin{align*}
W(\sqrt{x})&=\sum_{l\le \sqrt{x}} \frac{1}{l} (\sqrt{x}-l)^2 \sum_{q|l} qT(q)\\
&=\sum_{q\le\sqrt{x}}\sum_{r\le\sqrt{x}/q} \frac{1}{qr} (\sqrt{x}-qr)^2 qT(q)\\
&=\sum_{q\le\sqrt{x}} q^2T(q) \sum_{r\le\sqrt{x}/q} \frac{1}{r}\Big(\frac{\sqrt{x}}{q}-r\Big)^2\\
&=\sum_{q\le\sqrt{x}} q^2T(q) \Big(\frac{x}{q^2}\log\frac{\sqrt{x}}{q}+\frac{C_2x}{q^2}+\frac{\sqrt{x}}{q}+O(1)\Big)\\
&=\Big(\frac{1}{2}x\log x+C_2x\Big)\sum_{q\le\sqrt{x}} T(q)-x\sum_{q\le\sqrt{x}} T(q)\log q+\sqrt{x}\sum_{q\le\sqrt{x}} qT(q)+O(\sqrt{x}).
\end{align*}

For those sums, we have the following arguments which are similar to those in \S 5:

\begin{equation*}
\sum_{q>\sqrt{x}} T(q)\ll \frac{1}{\sqrt{x}} \text{, \quad}\sum_{q>\sqrt{x}} T(q)\log q\ll\frac{\log x}{\sqrt{x}} \text{\quad and \quad} \sum_{q\le\sqrt{x}} qT(q)\ll\log x.
\end{equation*}

So

\begin{equation*}
W(\sqrt{x})=\Big(\frac{1}{2}x\log x+C_2x\Big)\sum_{q=1}^{\infty} T(q)-x\sum_{q=1}^{\infty} T(q)\log q+O(\sqrt{x}\log x),
\end{equation*}

hence

\begin{equation}\label{xWofsqrtofxSpecial}
xW(\sqrt{x})=\frac{1}{2}\Big(\sum_{q=1}^{\infty} T(q)\Big) x^2\log x +\Big(C_2\sum_{q=1}^{\infty} T(q)-\sum_{q=1}^{\infty} T(q)\log q\Big) x^2+O(x^{3/2}\log x).
\end{equation}

Note that we have

\begin{equation}\label{constantD2}
D_2=\frac{1}{2}\Big(C_2\sum_{q=1}^{\infty} T(q)-\sum_{q=1}^{\infty} T(q)\log q\Big)
\end{equation}

when comparing (\ref{xWofsqrtofxSpecial}) with (\ref{xWofsqrtofx}).

\vspace{3mm}

By (\ref{2s6withwofXonly}), (\ref{vxqmaintermwiths6}), (\ref{constantc0}) and (\ref{constantc1}),

\begin{equation}\label{vofxqSpecial}
V(x,Q)=Q\sum_{n\le x}r_3(n)^2+x^2\Big(C_0\log\frac{x}{Q}+C_0C_2-C_1\Big)-\Gamma\Big(\frac{4}{3}\Big)^6Q^2 W\Big(\frac{x}{Q}\Big)+U(x,Q)+O(x^{3/2}\log x).
\end{equation}

According to the exponent of $x$, the error term $O(x^{3/2}\log x)$ can be absorbed into $U(x,Q)$. Finally, by definitions of the constants and the cases that were discussed at the beginning of this part, we have the following theorem and corollaries:

\vspace{3mm}


\textbf{Theorem 2.} \textit{Assume that $\varepsilon$ is a positive number and $k\ge1$ is an integer, and let}

\begin{equation*}
U_0(x,Q)=V(x,Q)-Q\sum_{n\le x} r_3(n)^2 - x^2\Big(C_0\log\frac{x}{Q}+C_0C_2-C_1\Big)-\Gamma\Big(\frac{4}{3}\Big)^6 \sum_{l\le k} \frac{h(l)}{l} (x-lQ)^2
\end{equation*}

\textit{where}

\begin{equation*}
C_0=\Gamma\Big(\frac{4}{3}\Big)^6 \sum_{q=1}^{\infty} \frac{1}{q^7}\sum_{\substack{c=1\\
(c,q)=1
}}^q |S(q,c)|^6,
\end{equation*}

\begin{equation*}
C_1=\Gamma\Big(\frac{4}{3}\Big)^6 \sum_{q=1}^{\infty} \frac{\log q}{q^7}\sum_{\substack{c=1\\
(c,q)=1
}}^q |S(q,c)|^6,
\end{equation*}

\begin{equation*}
C_2=-\frac{11}{12}-2\int^{\infty}_1 \frac{B_2(u)}{u^3} \mathrm{d}u,
\end{equation*}

\begin{equation*}
B_2(u)=\frac{1}{2}(u-[u])^2-\frac{1}{2}(u-[u])+\frac{1}{12},
\end{equation*}

\begin{equation*}
h(l)=\sum_{q|l}\frac{1}{q^6} \sum_{\substack{c=1\\
(c,q)=1
}}^q |S(q,c)|^6,
\end{equation*}

\textit{and}

\begin{equation*}
S(q,c)= \sum_{m=1}^q e\Big(\frac{cm^3}{q}\Big).
\end{equation*}

\textit{Then when $\displaystyle\frac{x}{k+1}<Q\le \frac{x}{k}$ one has}

\begin{equation*}
U_0(x,Q) \ll x^{\frac{10}{9}+\frac{2}{3}(A^+)+\varepsilon}
\end{equation*}

\textit{where $A^+$ is the upper exponent of $x$ in the representation of $\displaystyle\sum_{n\le x} r_3(n)^2$, or}

\begin{equation*}
U_0(x,Q) \ll x^{\frac{10}{9}+\varepsilon} \Big(\sum_{n\le x} r_3(n)^2\Big)^{\frac{2}{3}}
\end{equation*}

\textit{if the exponent of $x$ in the representation of $\displaystyle\sum_{n\le x} r_3(n)^2$ exists.}


\vspace{3mm}

\textbf{NOTE:} For the following corollaries, we use the same notations for constants and functions as in Theorem 2, and the error term $U_0(x,Q)$ or $U_0(x)$ has the same property as $U_0(x,Q)$ in Theorem 2.

\vspace{3mm}


\textbf{Corollary 1.} (i) \textit{When $\displaystyle\frac{x}{3}<Q\le\frac{x}{2}$ one has}

\begin{equation*}
V(x,Q)=Q\sum_{n\le x} r_3(n)^2 + x^2\Big(C_0\log\frac{x}{Q}+C_0C_2-C_1-\frac{3}{2}\Gamma\Big(\frac{4}{3}\Big)^6\Big)+4\Gamma\Big(\frac{4}{3}\Big)^6Qx-3\Gamma\Big(\frac{4}{3}\Big)^6 Q^2 +U_0(x,Q),
\end{equation*}

(ii) \textit{when $\displaystyle\frac{x}{2}<Q\le x$ one has}

\begin{equation*}
V(x,Q)=Q\sum_{n\le x} r_3(n)^2 + x^2\Big(C_0\log\frac{x}{Q}+C_0C_2-C_1-\Gamma\Big(\frac{4}{3}\Big)^6\Big)+2\Gamma\Big(\frac{4}{3}\Big)^6Qx-\Gamma\Big(\frac{4}{3}\Big)^6 Q^2 +U_0(x,Q).
\end{equation*}


\vspace{3mm}


\textbf{Corollary 2.} $\displaystyle V(x,x)=x\sum_{n\le x} r_3(n)^2 +x^2(C_0C_2-C_1)+U_0(x).$


\vspace{3mm}


\textbf{Corollary 3.} \textit{Assume that $m\ge1$ is a constant. Then one has}

\begin{equation*}
V\Big(x, \frac{x}{m}\Big)=\frac{x}{m}\sum_{n\le x} r_3(n)^2 +x^2\Big(C_0\log m+C_0C_2-C_1-\Gamma\Big(\frac{4}{3}\Big)^6\frac{1}{m^2}\sum_{l\le m} \frac{h(l)}{l} (m-l)^2 \Big) +U_0(x).
\end{equation*}


\vspace{3mm}

Note that the error terms in Theorem 2 and the corollaries are $\ll x^{17/9+\varepsilon}$ according to \S 8. That is why these results are more precise than Theorem 1.

\enddocument
\end{document}